\documentclass[12pt]{amsart}
\usepackage{amsmath,amssymb,amsthm,mathtools}
\usepackage[left=1in,right=1in,top=1.25in,bottom=1.25in]{geometry}
\usepackage[bookmarks=true,colorlinks,linkcolor=blue]{hyperref}
\usepackage{graphicx,subcaption}
\usepackage{mathrsfs}
\usepackage{standalone,tikz,tikz-cd}
\usetikzlibrary{positioning,knots,patterns,hobby,decorations.pathreplacing}
\usepackage{comment}
\allowdisplaybreaks
\usepackage{accents}
\usepackage{extarrows}

\newcommand{\bcon}{\begin{conjecture}}
\newcommand{\econ}{\end{conjecture}}
\newcommand{\bcor}{\begin{corollary}}
\newcommand{\ecor}{\end{corollary}}
\newcommand{\bdf}{\begin{definition}}
\newcommand{\edf}{\end{definition}}
\newcommand{\benu}{\begin{enumerate}}
\newcommand{\eenu}{\end{enumerate}}
\newcommand{\bexa}{\begin{example}}
\newcommand{\eexa}{\end{example}}
\newcommand{\bexe}{\begin{exercise}}
\newcommand{\eexe}{\end{exercise}}
\newcommand{\bfac}{\begin{fact}}
\newcommand{\efac}{\end{fact}}
\newcommand{\bite}{\begin{itemize}}
\newcommand{\eite}{\end{itemize}}
\newcommand{\blem}{\begin{lemma}}
\newcommand{\elem}{\end{lemma}}
\newcommand{\bmat}{\begin{pmatrix}}
\newcommand{\emat}{\end{pmatrix}}
\newcommand{\bprb}{\begin{problem}}
\newcommand{\eprb}{\end{problem}}
\newcommand{\bpro}{\begin{proposition}}
\newcommand{\epro}{\end{proposition}}

\newcommand{\bque}{\begin{question}}
\newcommand{\eque}{\end{question}}
\newcommand{\brem}{\begin{remark}}
\newcommand{\erem}{\end{remark}}
\newcommand{\bthm}{\begin{theorem}}
\newcommand{\ethm}{\end{theorem}}

\providecommand{\abs}[1]{\lvert #1\rvert}
\DeclareMathOperator{\rank}{rank}
\DeclareMathOperator{\tr}{tr}

\DeclareMathOperator{\Fr}{Fr}
\DeclareMathOperator{\Gr}{Gr}

\DeclareMathOperator{\id}{id}
\DeclareMathOperator{\pr}{pr}
\DeclareMathOperator{\dego}{\deg^\circ}
\DeclareMathOperator{\GKdim}{GK dim}
\newcommand{\qeq}{\stackrel{(q)}{=}}

\newcommand{\surface}{\mathfrak{S}}
\newcommand{\quasi}{\mathcal{E}}
\newcommand{\marked}{\mathcal{P}}
\newcommand{\face}{\mathcal{F}}

\makeatletter
\newcommand{\sbullet}{%
 \hbox{\fontfamily{lmr}\fontsize{.4\dimexpr(\f@size pt)}{0}\selectfont\textbullet}}
\DeclareRobustCommand{\mathbullet}{\accentset{\bullet}}
\makeatother

\newcommand{\fockY}{\mathcal{Y}}
\newcommand{\balancedY}{\mathcal{Y}^\mathrm{bl}}
\newcommand{\mullerX}{\mathfrak{X}}

\newcommand{\extYbl}{\bar{\mathcal{Y}}^\mathrm{bl}}
\newcommand{\extX}{\bar{\mathfrak{X}}}

\newcommand{\rd}{{\mathrm{rd}}}
\newcommand{\bad}{{\mathrm{bad}}}

\newcommand{\mullerS}{\mathscr{S}^+}
\newcommand{\stateS}{\mathscr{S}}
\newcommand{\reduceS}{\mathscr{S}^\rd}

\newcommand{\OSL}{\mathcal{O}_{q^2}(SL_2)}
\newcommand{\eqdef}{\overset{\mathrm{def}}{=\joinrel=}}

\theoremstyle{plain}
\newtheorem{theorem}{Theorem}[section]
\newtheorem{thm}{Theorem}
\newtheorem{lemma}[theorem]{Lemma}
\newtheorem{corollary}[theorem]{Corollary}
\newtheorem{proposition}[theorem]{Proposition}
\newtheorem{conjecture}{Conjecture}
\newtheorem{question}{Question}

\newtheorem{definition}{Definition}

\theoremstyle{definition}

\newtheorem{remark}[theorem]{Remark}
\newtheorem{example}[theorem]{Example}

\def\BC{\mathbb C}
\def\BN{\mathbb N}
\def\BZ{\mathbb Z}
\def\BR{\mathbb R}
\def\BT{\mathbb T}
\def\BQ{\mathbb Q}

\def\cR{{\mathcal R}}

\def\la{\langle}
\def\ra{\rangle}

\def\al{\alpha}
\def\ve{\varepsilon}

\def\be { \begin{equation} }
\def\ee { \end{equation} }

\def\bD{{\bar \Delta }}

\def\Dd{\Delta_\partial}

\def\tal{\tilde \al}

\def\bY{{\bar {\cY}}}
\def\bYbl{\bY^{\mathrm{bl}}}
\def\Ybl{\cY^{\mathrm{bl}}}

\newcommand\no[1]{}

\def\bP{\bar P}

\def\bT{\mathbb T}

\def\Do{{\mathring \Delta}}
\def\Pd{\cP_\partial}

\def\rk{\mathrm{rk}}

\def\cS{\mathscr S}
\def\ot{\otimes}

\def\cE{\mathcal E}
\def\cF{\mathcal F}

\def\bk{\mathbf k}

\def\bn{\mathbf n}

\def\cP{\mathcal P}

\def\fS{\mathfrak S}
\def\bfS{\overline{\fS}}

\def\cY{\mathcal Y}
\def\ev{{\mathrm{ev}}}

\def\embed{\hookrightarrow}

\def\sX{\mathfrak X}

\def\cY{\mathcal Y}
\def\bQ{\bar Q}
\def\ev{{\mathrm{ev}}}

\newcommand{\red}[1]{{\color{red}#1}}
\newcommand{\blue}[1]{{\color{blue}#1}}

\def\cB{{\mathfrak B}}

\def\cD{\mathcal D}
\def\CS{\cS(\fS)}
\def\onto{\twoheadrightarrow}
\def\embed{\hookrightarrow}
\def\pfS{\partial \fS}
\def\pbfS{\partial \bfS}
\def\cPd{\cP_\partial}
\def\cPo{\mathring{\cP}}

\def\ori{{\mathfrak o}}
\def\ofS{\mathring {\fS}}
\def\cSd{\cS^\rd}
\def\fM{\mathfrak{M}}
\def\fT{{\mathfrak{T}}}
\def\Po{\mathring {\cP}}
\def\bX{\bar {\sX}}
\def\bcE{\bar{\cE}}
\def\Ed{\cE_\partial}
\def\hEd{\widehat{\Ed}}

\begin{document}

\title{Quantum traces and embeddings of stated skein algebras into quantum tori}
\author[Thang T. Q. L\^e]{Thang T. Q. L\^e}
\address{School of Mathematics, 686 Cherry Street,
 Georgia Tech, Atlanta, GA 30332, USA}
\email{letu@math.gatech.edu}
\author[Tao Yu]{Tao Yu}
\address{School of Mathematics, 686 Cherry Street,
 Georgia Tech, Atlanta, GA 30332, USA}
\email{tyu70@gatech.edu}


\thanks{
2010 \emph{Mathematics Classification:} Primary 57N10. Secondary 57M25.\\
}

\begin{abstract}
The stated skein algebra of a punctured bordered surface (or equivalently, a marked surface) is a generalization of the well-known Kauffman bracket skein algebra of unmarked surfaces and can be considered as an extension of the quantum special linear group $\OSL$ from a bigon to general surfaces.

We show that the stated skein algebra of a punctured bordered surface with non-empty boundary can be embedded into quantum tori in two different ways.
The first embedding can be considered as a quantization of the map expressing the trace of a closed curve in terms of the shear coordinates of the enhanced Teichm\"uller space, and is a lift of Bonahon-Wong's quantum trace map.
The second embedding can be considered as a quantization of the map expresses the trace of a closed curve in terms of the lambda length coordinates of the decorated Teichm\"uller space, and is an extension of Muller's quantum trace map. We explain the relation between the two quantum trace maps. We also show that the quantum cluster algebra of Muller is equal to a reduced version of the stated skein algebra.
As applications we show that the stated skein algebra is an orderly finitely generated Noetherian domain and calculate its Gelfand-Kirillov dimension. 
\end{abstract}
\maketitle

\tableofcontents

\def\SS{\cS(\fS)}
\def\SSo{\mathring{\cS}(\fS)}

\section{Introduction}
\subsection{Kauffman bracket skein algebra and Bonahon-Wong's quantum trace}
Let $\fS$ be an oriented surface.
The ground ring $\cR$ is a commutative Noetherian domain with 1, containing a distinguished invertible element $q^{1/2}$. For example $\cR=\BZ[q^{\pm 1/2}]$, or $\cR=\BC$ and $q^{1/2}$ is a non-zero complex number.

The Kauffman bracket skein algebra $\SSo$, introduced by Przytycki \cite{Przy} and Turaev \cite{Turaev,Turaev2}, is the $\cR$-module freely generated by framed unoriented links in the thickened surface $
\fS \times (-1,1)$ subject to the Kauffman skein relations
\begin{align*}
\begin{tikzpicture}[scale=0.8,baseline=0.3cm]
\fill[gray!20!white] (-0.1,0)rectangle(1.1,1);
\begin{knot}[clip width=8,background color=gray!20!white]
\strand[very thick] (1,1)--(0,0);
\strand[very thick] (0,1)--(1,0);
\end{knot}
\end{tikzpicture}
&=q
\begin{tikzpicture}[scale=0.8,baseline=0.3cm]
\fill[gray!20!white] (-0.1,0)rectangle(1.1,1);
\draw[very thick] (0,0)..controls (0.5,0.5)..(0,1);
\draw[very thick] (1,0)..controls (0.5,0.5)..(1,1);
\end{tikzpicture}
+q^{-1}
\begin{tikzpicture}[scale=0.8,baseline=0.3cm]
\fill[gray!20!white] (-0.1,0)rectangle(1.1,1);
\draw[very thick] (0,0)..controls (0.5,0.5)..(1,0);
\draw[very thick] (0,1)..controls (0.5,0.5)..(1,1);
\end{tikzpicture}\, ,\qquad
\begin{tikzpicture}[scale=0.8,baseline=0.3cm]
\fill[gray!20!white] (0,0)rectangle(1,1);
\draw[very thick] (0.5,0.5)circle(0.3);
\end{tikzpicture}
=(-q^2 -q^{-2})
\begin{tikzpicture}[scale=0.8,baseline=0.3cm]
\fill[gray!20!white] (0,0)rectangle(1,1);
\end{tikzpicture}\, .
\end{align*}
The product is given by stacking. See Section \ref{sec.surfaces} for details.

The skein algebra $\SSo$ and its analogs have played an important role in low dimensional topology and quantum algebra as they have applications and connections to objects such as character varieties \cite{Bullock,PS,Turaev,BFK,CM}, the Jones polynomial and its related topological quantum field theory (TQFT) \cite{Kauffman,BHMV,Tu:Book}, (quantum) Teichm\"uller spaces and (quantum) cluster algebras \cite{BW1,FGo, Muller}, the AJ conjecture \cite{FGL,Le5}, and many more. It is important to understand algebraic properties of the skein algebra $\SSo$, for example its representation theory.

In an important development, Bonahon and Wong \cite{BW1} proved that when $\fS$ is the result of removing at least one puncture from a closed surface, the skein algebra $\SSo$ can be embedded into a \emph{quantum torus} by the \emph{quantum trace map}
\[\tr: \SSo \embed \bT(Q)\]
which is a quantization of the map expressing the trace of curves in the shear coordinates of Teichm\"uller space. Here the quantum torus of an antisymmetric $r\times r$ integral matrix $Q$ is the $\cR$-algebra
\[\bT(Q):= \cR\la x_i^{\pm 1}, i=1, \dots, r\ra /(x_i x_j = q^{Q_{ij}} x_j x_i).\]
Thus $\bT(Q)$ is the algebra of Laurent polynomials in the $r$ variables $x_i$ which might not commute but are $q$-commuting in the sense that $x_i x_j= q^{Q_{ij}} x_j x_i$. The quantum torus $\bT(Q)$ is known to a be a Noetherian domain with Gelfand-Kirillov dimension $r$ and is well-understood in algebraic terms. The subalgebra generated by non-negative powers of $x_i$ is called the {\em quantum space}. As a Noetherian domain, $\bT(Q)$ has a ring of fractions $\Fr(\bT(Q))$ which is a division algebra.

To every ideal triangulation $\Delta$ of the punctured surface $\fS$, one can associate a quantum torus $\Ybl(\fS;\Delta)$, which is the square root version of the Chekhov-Fock algebra \cite{BW1,Hiatt}, known as the quantum Teichm\"uller space \cite{CF,Kashaev2}, as it quantizes the enhanced Teichm\"uller space of the surface. The Bonahon-Wong quantum trace embeds the skein algebra $\SSo$ into the quantum torus $\Ybl(\fS;\Delta)$. 
For another triangulation there is a coordinate change isomorphism between the two rings of fractions of the Chekhov-Fock algebras which permutes the quantum trace maps.

The quantum trace helps to understand the skein algebra algebraically, and opens possibilities to quantize Thurston's theory of hyperbolic surfaces to build hyperbolic topological field theory and to better understand the volume conjecture \cite{Kashaev1,MM}.


\subsection{Stated skein algebra}

To better understand the quantum trace map, the first author~\cite{Le:TDEC} introduced the {\em stated skein algebra} $\SS$ for a {\em punctured bordered surface $\fS$}, which is the result of removing a finite number of punctures from a compact oriented surface. We assume that every connected component of the boundary $\pfS$ is diffeomorphic to the open interval $(0,1)$. The stated skein algebra $\SS$ is a quotient of a coarser version considered by Bonahon and Wong. In addition to framed links, framed tangles properly embedded into the thicken surface $\fS \times (-1,1)$ are allowed, and new relations are introduced at the boundary. The tangles are equipped with states, which assign a positive sign or a negative sign to each endpoint of the tangles. For details, see Section \ref{sec.surfaces}.

A main feature of the stated skein algebra is the existence of the splitting homomorphism
\begin{equation}
\theta_c: \SS \to \cS(\fS'),
\end{equation}
where $\fS'$ is the result of slitting $\fS$ along an ideal arc $c$. The splitting map $\theta_c$ is an algebra embedding and is given by a simple state sum. The precise image is described by a certain invariant subspace in terms of Hochschild cohomology \cite{CL,KQ}. The splitting map gives a new, simple proof of the existence of the quantum trace map.

\def\fB{{\mathfrak B}}
The stated skein algebra theory has a rich mathematical content, far beyond the application in the understanding of the quantum trace, as it fits well with the \emph{integral} quantum group associated to $SL_2(\BC)$ and its integral dual $\mathcal O_{q^2}(SL_2)$. Many algebraic facts concerning the quantum groups have simple transparent interpretations by geometric formulas via the theory of stated skein algebra. For example, the stated skein algebra $\cS(\fB)$ of the bigon $\fB$, with its natural cobraided structure, is isomorphic to the cobraided Hopf algebra $\OSL$, and under the isomorphism the natural basis of $\cS(\fB)$ maps to Kashiwara's canonical basis \cite{Ka} of $\OSL$, see \cite{CL}. The stated skein algebra has connections to the moduli algebra of Alekseev-Gross-Schomerus \cite{AGS} and Buffenoir-Roche \cite{BR} and factorization homology of surfaces \cite{BBJ}, see \cite{Faitg, Korinman, LY0}, and was further studied in \cite{CL,KQ,BL}, among others. There are generalizations to other Lie algebras \cite{Higgins,LS,LY1}.

Suppose $\Delta$ is an ideal triangulation of $\fS$. One can define the Chekhov-Fock algebra $\Ybl(\fS;\Delta)$. Bonahon-Wong's quantum trace, defined on a coarser version of the stated skein algebra, actually gives an algebra homomorphism
\begin{equation}\label{eq.22}
\tr_\Delta: \SS \to \Ybl(\fS;\Delta),
\end{equation}
which is only injective when the boundary $\pfS$ is empty. As embeddings into quantum tori tell us a lot about the algebraic structure of $\SS$, we are looking for such embeddings for surfaces with boundaries.


\def\SSp{\cS^+(\fS)}
\def\bal{{\mathrm{bl}}}
\def\bYbl{\bY^\bal}
\def\bTheta{\bar \Theta}
\def\cZ{\mathcal Z}

\def\vpE{\varphi_\cE}
\def\bXE{\bX(\fS;\cE)}
\def\bpXE{\bX^\diamond (\fS;\cE)}
\def\bYD{\bYbl(\fS;\Delta)}
\def\pD{\phi_\Delta}
\def\bE{\bar {\cE}}
\def\RVb{\cR[\Po]^\diamond }
\def\vpEb{\varphi_\cE^\diamond }
\def\bXEb{\bX^\diamond (\fS;\cE)}
\def\CSb{\cS^\diamond (\fS)}

\def\YD{\bYbl(\fS;\Delta)}

\subsection{Shear coordinate quantum trace map}

Assume $\Delta$ is an ideal triangulation of a punctured bordered surface $\fS$. We will define a bigger quantum torus $\YD$ that comes with an $\cR$-linear projection $\pr: \YD \onto \Ybl(\fS;\Delta)$ onto the original Chekhov-Fock algebra used in Bonahon-Wong's quantum trace map \eqref{eq.22}. 

\begin{thm}[part of Theorem \ref{thm.embed1}]\label{thm.1aa}
There is a natural algebra embedding
\begin{equation}\label{eq.qtr2a}
\phi_\Delta: \SS \embed \YD
\end{equation}
which lifts the Bonahon-Wong quantum trace map \eqref{eq.22}, i.e., one has $\tr_\Delta= \pr\circ \phi_\Delta$.
\end{thm}

Here, naturality is with respect to triangulation changes. This means for another triangulation $\Delta'$, there is a coordinate change isomorphism $\bTheta_{\Delta'\Delta}: \Fr(\YD) \to \Fr(\bYbl(\fS;\Delta'))$ which transfers $\phi_\Delta$ to $\phi_{\Delta'}$. Here $\Fr(A)$ is the ring of fractions of the Ore domain $A$. The coordinate change isomorphism is functorial, and it extends the original coordinate change isomorphism of the Chekhov-Fock algebra defined in \cite{Hiatt,BW1,CF,Liu}.

Using the top degree part of the extended quantum trace $\phi_\Delta$, we will establish the following facts about the stated skein algebra of a punctured bordered surface.

\begin{thm}[See Theorem \ref{thm-noether1}]
\label{thm-noether10}
Let $\fS$ be a punctured bordered surface. 
\begin{enumerate}
\item The $\cR$-algebra $\cS(\fS)$ is orderly finitely generated. This means, there are elements $\alpha_1, \dots, \alpha_n\in \SS$ such that the set $\{ \alpha_1^ {k_1} \dots \alpha_n^{k_n} \mid k_i \in \BN\}$ spans $\cS(\fS)$ over $\cR$.
\item $\cS(\fS)$ is a Noetherian domain.
\item If $\fS$ has an ideal triangulation, then
the Gelfand-Kirillov dimension of $\cS(\fS)$ is
\begin{equation}\label{eq.rs}
r(\fS)\eqdef 3 |\Pd| - 3 \chi(\fS),
\end{equation}
where $\Pd$ is the set of boundary punctures and $\chi(\fS)$ is the Euler characteristic.
\end{enumerate}
\end{thm}


Part (c) of Theorem~\ref{thm-noether10} will imply that the domain and the target space of the embedding $\pD$ in \eqref{eq.qtr2a} have the same Gelfand-Kirillov dimension. Hence the embedding $\pD$ is tight in the sense that its image cannot lie in a quantum torus of lower dimension.


If $\pfS=\emptyset$ then $\SS$ is the ordinary skein algebra, and most of Theorem \ref{thm-noether1} was known: the finite generation (without order) was proved in \cite{Bullock}, the orderly finite generation was proved in \cite{AF}, and the Noetherian domain property was established in \cite{PS2}.

\subsection{Length coordinate quantum trace map}

Suppose the punctured bordered surface $\fS$ has non-empty boundary. Except for a few simple cases, $\fS$ has a {\em quasitriangulation}, which is a collection of ideal arcs whose endpoints are on the boundary and which cuts $\fS$ into ideal triangles and once punctured monogons, see Section \ref{sec.length}. We will associate to a quasitriangulation $\cE$ a quantum torus $\bX(\fS;\cE)$ which contains the quantum space $\bX_+(\fS;\cE)$.
\begin{thm}[part of Theorem \ref{thm.embed1}]
There is a natural algebra embedding
\[\varphi_\cE: \SS \embed \bX(\fS;\cE)\]
such that the image $\varphi_\cE(\SS)$ is sandwiched between $\bX_+(\fS;\cE)$ and $\bX(\fS;\cE)$. The algebra $\SS$ is an Ore domain and its ring of fractions $\Fr(\SS)$ is a division ring isomorphic to the ring of fractions $\Fr(\bX(\fS;\cE))$.
\end{thm}
The embedding map $\vpE$ has a simple geometric interpretation.
As $\SS$ is sandwiched between the quantum space and the quantum torus, one can derive many consequences. For example, the center of $\SS$ is the intersection of $\SS$ and the center of the quantum torus $\sX(\fS;\cE)$, which is easy to calculate. Representations of $\SS$ can be studied from this point of view. In this aspect the embedding $\vpE$ is better than the embedding $\pD$ of Theorem~\ref{thm.1aa}.

It was observed in \cite{Le:TDEC} that the subspace $\SSp$ of $\SS$ spanned by tangles with positive states is isomorphic to the Muller skein algebra \cite{Muller}. For the Muller subalgebra $\SSp$, the map $\varphi_\cE$, with a smaller target space, was constructed in \cite{Muller} in the case when there is no interior puncture, and in \cite{LP} for the general case. When there is no interior puncture and $q=1$, the space $\SSp$ can be identified with a space of functions on the decorated Teichm\"uller space, and $\varphi_\cE(\al)$ for a loop $\al$ expresses the trace of $\al$ in the holonomy representation of the hyperbolic metric in terms of the lambda lengths of the edges of the triangulation. Thus our embedding $\varphi_\cE$ is an extension/generalization of Muller's quantum trace map to the full stated skein algebra. It should be noted that the additional negative states present non-trivial obstacles that we have to overcome.

\subsection{Relation between the two quantum trace map}
Suppose $\fS$ has non-empty boundary, and $\cE$ is a quasitriangulation of $\fS$. There is a unique ideal triangulation $\Delta$ which is an extension of $\cE$. We now have two quantum trace maps, the algebra embeddings $\pD$ and $\vpE$.

\def\psE{\psi_\cE}
\begin{thm}[part of Theorem \ref{thm-len2shear}]\label{thm.equi2}
There is an $\cR$-algebra embedding
\[\psi_\quasi:\bXE\xhookrightarrow{ } \extYbl(\surface;\Delta)\]
 making the following diagram commute:
\[\begin{tikzcd}[row sep=tiny]
&\extYbl(\surface;\Delta)\\
\stateS(\surface)\arrow[ru,"\phi_\Delta"]
\arrow[rd,"\varphi_\quasi"']&\\
&\extX(\surface;\quasi)\arrow[uu,"\psi_\quasi"]
\end{tikzcd}\]
\end{thm}
The map $\psE$ has a transparent, simple geometric interpretation.
Moreover, we show that $\extYbl(\surface;\Delta)$ is a central quadratic extension of $\extX(\surface;\quasi)$. When $\fS$ has no interior puncture and $q=1$, the map $\psE$ is essentially the well-known map which changes shear coordinates to lambda-length coordinates in Teichm\"uller spaces.

\def\cA{\mathcal A}
\subsection{Quantum cluster algebra of Muller}
For the case when $\fS$ has no interior puncture and non-empty boundary, Muller \cite{Muller} constructed a quantum cluster algebra $\cA(\fS)$, which is a quantization of the cluster algebra associated to $\fS$ defined in \cite{GSV,FG}. The quantum cluster algebra $\cA(\fS)$ is a subalgebra of the localization $\SSp\fM^{-1}$, where $\fM$ is the multiplicative subset generated by boundary edges. Recall that $\SSp$, the Muller skein algebra, is the submodule of $\SS$ generated by tangles with positive states only. When $\fS$ has at least two boundary punctures, Muller showed that $\SSp\fM^{-1}=\cA(\fS)$.

The kernel of Bonahon-Wong's quantum trace \eqref{eq.22} is calculated in \cite{CL}, and is generated by the so called bad arcs, see Section \ref{sec.length}. By factoring out the kernel, one gets the {\em reduced skein algebra} and an algebra embedding
\begin{equation}
\tr^\rd_\Delta: \cS^\rd (\fS) \embed \Ybl(\fS;\Delta).
\end{equation}

Consider the composition
\[\varkappa: \SSp \embed \SS \onto \cS^\rd(\fS)\]
which maps the Muller algebra $\SSp$ to the reduced skein algebra.

\def\tk{\tilde \varkappa}

\begin{thm}[See Theorem \ref{thm.reduced-cluster}]\label{thm.rdclus}
Suppose $\fS$ is a connected punctured bordered surface. The embedding $\varkappa$ extends uniquely to an $\cR$-algebra isomorphism
$\tk: \mullerS(\surface)\fM^{-1} \xrightarrow{\cong} \reduceS(\surface)$.

Consequently, when $\fS$ has at least two punctures but no interior puncture, the quantum cluster algebra $\cA(\fS)$ is naturally isomorphic to the reduced skein algebra $\reduceS(\surface)$.
\end{thm}

While $\mullerS(\surface)$ and its localization $\mullerS(\surface)\mathfrak{M}^{-1}$ are defined using only positive states, the reduced stated skein algebra $\reduceS(\surface)$ uses both positive and negative states. Hence it is a surprise that we can have the above result, which demonstrates the ubiquity of the stated skein algebra and gives a new perspective for the quantum cluster algebra of surfaces.

\subsection{Applications, related works}
Let $\cR=\BC$ and $q$ be a root of 1. For applications in hyperbolic TQFT theory, one would like to know representations theory of $\SS$. For this one needs to know the center of $\SS$, and dimension of $\SS$ over its center. Besides $\SS$ is a Poisson order \cite{BG}, and to understand the representations of $\SS$ one would like to calculate the Poisson structure on the center of $\SS$. All these problems can be approached using the embeddings of $\SS$ into quantum tori, and we will explore these questions in upcoming work.

Our embedding of the stated skein algebra into quantum tori is close to, and actually related to, results about embedding of quantized enveloping algebras into quantum tori \cite{Faddeev,SS}. Some geometric ideas and techniques in the current paper can be applied to higher ranked Lie algebra in our upcoming work \cite{LY1} concerning the quantum trace of $SL_n$ skein algebra.

\subsection{Acknowledgments}
The authors would like to thank F. Bonahon, F. Costantino, C. Frohman, J. Kania-Bartoszynska, A. Kricker, A. Sikora, and M. Yakimov for helpful discussions. The first author is supported in part by NSF grant DMS 1811114, and benefited from a visit to Nangyang Technological University in November 2019.

\subsection{Organization of the paper}
In Section \ref{sec.alg} we fix some notations and explain some algebraic facts. Section \ref{sec.surfaces} contains basics of stated skein algebras. Section \ref{sec.length} has the embedding of $\SS$ into a quantum torus which is the quantization of length coordinate functions. Section \ref{sec.red} discusses the reduced skein algebra and the quantum cluster algebra. Section \ref{sec-shear} presents the embedding of $\SS$ into a quantum torus which is the quantization of shear coordinate functions. Section \ref{sec.rel} explains the relation between the two quantum trace maps.

\section{Notations, algebraic preliminaries}\label{sec.alg}

We fix notations and review the theory of quantum tori and Gelfand-Kirillov dimension.

\subsection{Notations, conventions}
Throughout the paper the ground ring $\cR$ is a commutative Noetherian domain with unit, with a distinguished invertible element $q^{1/2}$. All algebras are $\cR$-algebras unless otherwise stated.

Two elements $x,y$ in an $\cR$-algebra $A$ are \emph{$q$-proportional}, denoted by $x\qeq y$, if there is $k\in \BZ$ such that $x = q^{k/2} y$. Two elements $x,y\in A$ are \emph{$q$-commuting} if $xy$ and $yx$ are $q$-proportional

We denote by $\BN, \BZ, \BC$ respectively the set of non-negative integers, the set of integers, and the set of complex numbers. We emphasize that our $\BN$ contains 0.

\subsection{Weyl ordering}
 Suppose $x_1,x_2,\dots,x_n$ are pairwise $q$-commuting elements, $x_i x_j = q^{c_{ij}} x_j x_i$.
 The well-known \emph{Weyl normalization} of the product $x_1x_2\dots x_n$ is
\[[x_1x_2\dots x_n]=q^{-\frac{1}{2}\sum_{i<j}c(x_i,x_j)}x_1x_2\dots x_n.\]
The factor is chosen such that the normalization does not depend on the order of the product. That is, if $\sigma$ is a permutation of $\{1,2,\dots,n\}$, then $[x_1x_2\dots x_n]=[x_{\sigma(1)}x_{\sigma(2)}\dots x_{\sigma(n)}]$.

\subsection{Quantum torus}
The \emph{quantum torus} associated to an antisymmetric $r\times r$ integral matrix $Q$ is the algebra
\[\mathbb{T}(Q)\eqdef \cR\langle x_1^{\pm1},\dots,x_r^{\pm1}\rangle/\langle x_ix_j=q^{Q_{ij}}x_jx_i\rangle.\]

For $\bk=(k_1,\dots, k_r)\in \BZ^r$, let
\[ x^\bk \eqdef [ x_1 ^{k_1} x_2^{k_2} \dots x_r ^{k_r} ]= q^{-\frac{1}{2} \sum_{i<j} Q_{ij} k_i k_j} x_1 ^{k_1} x_2^{k_2} \dots x_r ^{k_r}\]
be the (Weyl) normalized monomial. Then $\{ x^\bk \mid \bk \in \BZ^r\}$ is a free $\cR$-basis of $\bT(Q)$, and
\begin{align}
\label{eq.prod}
x^\bk x ^{\bk'}& = q ^{\frac 12 \la \bk, \bk'\ra_Q} x^{\bk + \bk'},\quad \text{where } \la \bk, \bk'\ra_Q := \sum_{1\le i, j \le r} Q_{ij} k_i k'_j,\\
\label{eq.commu}
x^\bk x ^{\bk'} & = q ^{ \la \bk, \bk'\ra_Q} x^{\bk '} x ^\bk.
\end{align}
It follows that the decomposition
\begin{equation}\label{eq.grad}
\bT(Q) = \bigoplus_{\bk\in \BZ^r} \cR x^\bk
\end{equation}
gives the algebra $\bT(Q)$ a $\BZ^r$-grading.
A quantum torus is a Noetherian domain \cite{GW}. In particular, it has a ring of fractions, denoted by $\Fr(\bT(Q))$, which is a division algebra.

Suppose $Q'$ is another antisymmetric $r'\times r'$ integral matrix such that $HQ' H^T= Q$, where $H$ is an $r\times r'$ integral matrix and $H^T$ is its transpose. Then the $\cR$-linear map $\bT(Q)\to \bT(Q')$ given on the basis by $x^\bk \mapsto x^{\bk H}$ is an algebra homomorphism, called a \emph{multiplicatively linear homomorphism}. Here $\bk H$ is the product of the row vector $\bk$ and the matrix $H$.

When $\cR=\BZ[q^{\pm 1/2}]$, there is a unique $\BZ$-linear ring anti-automorphism
\[\omega:\mathbb{T}(Q)\to\mathbb{T}(Q), \quad \text{given by} \quad
\omega(q^{1/2})=q^{-1/2},\quad \omega(x_i)=x_i.\]
Here anti-automorphism means isomorphism to the opposite ring. Thus, $\omega(x+y)=\omega(x)+\omega(y)$, $\omega(xy)=\omega(y)\omega(x)$. Clearly $\omega^2=\operatorname{id}$. We will call $\omega$ the \emph{reflection anti-involution} of the quantum torus. An element $z\in \bT(Q)$ is \emph{reflection invariant} if $\omega(z)=z$. A map $f:\mathbb{T}(Q)\to\mathbb{T}(Q')$ is \emph{reflection invariant} if $f\circ\omega=\omega\circ f$. For example, all normalized monomials $x^\bk$ are reflection invariant, and all multiplicatively linear homomorphisms are reflection invariant.

\subsection{Monomial subalgebra}
Still assume $Q$ is an antisymmetric integral $r\times r$ matrix.
If $\Lambda\subset \BZ^r$ is a submonoid, then the $\cR$-submodule
$A(Q;\Lambda)\subset \bT(Q)$ spanned by $\{ x^\bk\mid \bk \in \Lambda \}$ is an $\cR$-subalgebra of $\bT(Q)$, called a \emph{monomial subalgebra}. When $\Lambda= \BN^r$, the corresponding subalgebra is the \emph{quantum space}, denoted by $\bT_+(Q)$, which is Noetherian \cite{BG}.

\begin{lemma}\label{r.mono}
If $\Lambda\subset \BZ^r$ is a submonoid finitely generated as an $\BN$-module, then the monomial subalgebra $A(Q;\Lambda)\subset \bT(Q)$ is a Noetherian domain.
\end{lemma}

\begin{proof}
As $A(Q;\Lambda)$ is a subalgebra of $\bT(\BQ)$, it is a domain. Suppose $\{\bn_1,\dots,\bn_k\}$ is an $\BN$-spanning set of $\Lambda$. Then $\Lambda$ is the quotient of the free $\BN$-module with basis $\{\bn_1,\dots, \bn_k\}$. Let $P_{ij}= \la \bn_i, \bn_j\ra_{Q}$. Then $A(Q;\Lambda) $ is a quotient of $\BT_+(P)$, and hence Noetherian.
\end{proof}

Note when $\Lambda$ is a subgroup, then $A(Q;\Lambda)$ is a quantum torus. This occurs as the image of a multiplicatively linear homomorphism.

Another example of a monomial subalgebra is the center.
As $\bT(Q)$ is graded by \eqref{eq.grad}, its center $Z(\bT(Q))$ is a graded subalgebra. In other words, there is a subgroup $\Lambda(Q,q)\subset \BZ^r$ such that $Z(\bT(Q))$ is the $\cR$-linear span of $\{x^\bk \mid \bk \in \Lambda(Q,q)\}$. That is $Z(\bT(Q)) = A(Q;\Lambda(Q,q))$. From the commutation relation \eqref{eq.commu} we see that
\[\Lambda(Q,q)= \{ \bk \in \BZ^r \mid q^{\la \bk, \bk' \ra_Q} =1\text{ for all } \bk' \in \BZ^r\}.\]

\subsection{Embedding into quantum torus}

\begin{proposition}\label{r.Ore}
Let $Q$ be an antisymmetric integral $r\times r$ matrix, and let $A$ be an $\cR$-algebra containing the quantum space $\bT_+(Q)$ as a subalgebra. Assume that $A$ is a domain, and for every $a\in A$ there is $\bk\in \BN^r$ such that $x^\bk a\in \bT_+(Q)\subset A$.

Then $A$ is an Ore domain, and the embedding $\bT_+(Q) \embed \bT(Q)$ can be uniquely extended to an algebra embedding $A \embed \bT(Q)$ which induces an isomorphism of the rings of fractions $\Fr(A) \xrightarrow{\cong} \Fr(\bT(Q))$.
\end{proposition}

\begin{proof}
Uniqueness is obvious. From \eqref{eq.prod} we see that for every $\bk\in \BN^r$ the map
\[\tau_\bk: \BT(Q)\to \BT(Q),\quad \tau_\bk(a)= x^{\bk} a x^{-\bk}\]
is an algebra automorphism which preserves $\BT_+(Q)$. If $x^\bk a =u \in \BT_+(Q)$, then 
\[ x^\bk ( a x^{\bk} - (\tau_\bk)^{-1}(u))= x^\bk u - x^\bk u =0.\]
Since $A$ is a domain, we have $a x^{\bk} = (\tau_\bk)^{-1}(u)\in \BT_+(Q)$.

We now show that the multiplicative subset $S= \{ q^{n/2} x^\bk \mid n\in \BZ, \bk\in \BN^r\}$ is an Ore set of $A$.
Let $a\in A$ and $\bk\in \BN^r$. To show that $S$ is right Ore, we have to show that $ A x^\bk \cap Sa\neq \emptyset$. By assumption $x^\bn a \in \BT_+(Q)$ for some $\bn\in \BN^r$. The following element
\[(x^\bk x^\bn) a = x^\bk (x^\bn a) = \tau_\bk (x^\bn a) x^{\bk}\]
belongs to both $ Sa$ and $A x^\bk $. Thus $S$ is a right Ore multiplicative subset of $A$. Similarly $S$ is also a left Ore subset. The right Ore localization $A[S^{-1}]$ contains $A$ as a subset since $S$, as a subset of $A$, does not have zero divisor.

It is clear that $S$ is also a right Ore set of $\BT_+(Q)$. Since localization is exact, we have an embedding
\begin{equation}
\gamma: \BT(Q)=\BT_+(Q)[S^{-1}] \embed A[S^{-1}].
\end{equation}
Since for every $a\in A$ there is $\bk\in \BN^r$ such that $a x^\bk\in \BT_+(Q)$, we see $\gamma$ is surjective, and hence $\gamma$ is an isomorphism. The restriction of $\gamma^{-1}$ onto $A$ is the desired extension.

As $\Fr(\bT_+(Q))= \Fr(\bT(Q))$ and $A$ is sandwiched between $\bT_+(Q)$ and $\bT(Q)$, we see that $A$ is an Ore domain, and $\Fr(A)=\Fr(\bT(Q))$.
\end{proof}

\subsection{Gelfand-Kirillov dimension}
The Gelfand-Kirillov dimension is a noncommutative analog of the Krull dimension.
Let $A$ be a finitely generated algebra over a field $k$, and let $V$ be a finite dimensional generating subspace, e.g., the span of the generators. Set $V^0=k$, $V^n=\{a_1 a_2\dots a_n\mid a_i\in V,i=1,\dots,n\}$ and $A_n=\sum_{i=0}^n V^i$. The \emph{Gelfand-Kirillov dimension}, or GK dimension,
is defined as
\[\GKdim A \eqdef \limsup_{n\to\infty}\frac{\log \dim_k (A_n)}{\log n}.\]
The dimension is independent of the choice of $V$. We extend the definition to an $\mathcal{R}$-algebra $A$ by letting $k=\Fr(\mathcal{R})$, the field of fractions of $\mathcal{R}$, and
\[\GKdim A \eqdef \GKdim(A\otimes_{\mathcal{R}} k).\]

\begin{lemma}\label{lemma-GKdim}
Let $A$ be a finitely generated $\cR$-algebra.
\begin{enumerate}
\item If $B$ is a (finitely generated) subalgebra or a quotient of $A$, then $\GKdim B\le\GKdim A$.
\item Suppose $\{F_k\}_{k=0}^\infty$ is a finite dimensional filtration of $A$, then the associated graded algebra $\Gr A$ has the same dimension $\GKdim(\Gr A)=\GKdim A$.
\item The quantum torus $\mathbb{T}(Q)$ and the quantum space $\mathbb{T}_+(Q)$ has GK dimension $r$, the number of generators.
\item More generally, the GK dimension of the monomial subalgebra $A(Q;\Lambda)$ is $\rank\Lambda$.
\end{enumerate}
\end{lemma}

\begin{proof}
(a) and (b) can be found in \cite[Propositions 8.2.2, 8.6.5]{MR}. (c) can be found in \cite[Theorem 25]{Reyes}.

For (d), let $\bar{\Lambda}\subset\mathbb{Z}^r$ be the subgroup generated by $\Lambda$. Thus $A(Q;\Lambda)$ is contained in the quantum torus $A(Q;\bar{\Lambda})$ with $k=\rank\Lambda$ generators. $\Lambda$ contains a free $\BN$-module of rank $k$. Hence the monomial algebra $A(Q, \Lambda)$ contains an isomorphic image of a quantum space with $k$ generators. By (a) and (c), all three algebras have GK dimension $k$.
\end{proof}

\def\ba{\bar a}
\section{Stated skein algebra}\label{sec.surfaces}

\subsection{Punctured bordered surfaces}
A \emph{punctured bordered surface} $\fS$ is a surface of the form $\fS = \bfS \setminus \cP$, where $\bfS$ is a compact oriented surface with
(possibly empty) boundary $\pbfS$, and $\cP$ is a finite set such that every connected component of the boundary $\pbfS$ has at least one point in $\cP$. We do not require $\fS $ be to connected. Each connected component of the boundary $\pfS$ is an open interval called a \emph{boundary edge} of $\fS$.

Points in $\cP$ are called \emph{punctures} of $\fS$. A puncture belonging to the boundary $\pbfS$ is called a \emph{boundary puncture}, and the rest are \emph{interior punctures}. Let $\cPd$ and $\cPo$ be the set of all boundary punctures and the set of all interior punctures respectively.

An \emph{ideal $n$-gon} is the standard disk with $n$ points on its boundary removed. For $n=1,2,3,4$, we also call them \emph{monogon}, \emph{bigon}, \emph{triangle} and \emph{quadrilateral}, respectively.

An \emph{ideal arc} of $\fS$ is an embedding $a: (0,1) \to \fS$ that can be extended to an immersion $\ba: [0,1]\to \bfS$ with $\ba(0),\ba(1)\in \cP$. We call $\ba(0)$ and $ \ba(1)$, which might coincide, the \emph{ideal endpoints} of $a$. As usual we identify an ideal arc with its image, which is considered as a non-oriented 1-dimensional submanifold. Isotopies of ideal arcs are considered in the class of ideal arcs.
If $\ba(0)=\ba(1)$ and $\ba$ bounds a disk in $\fS$, $a$ is called a \emph{trivial ideal arc}.

\def\tpfS{\partial \widetilde{\fS}}
\def\tfS{\widetilde{\fS}}
\subsection{Tangles in thickening of surface}

Fix a punctured bordered surface $\fS$. The \emph{thickening} of $\fS$ is the 3-manifold $\tfS= \fS \times (-1,1)$ with the orientation induced from those of $\fS$ and $(-1,1)$. If $b$ is a boundary edge of $\fS$ then $b \times (-1,1)$ is called a \emph{boundary wall} of $\tfS$. The union of all the boundary walls is the boundary $\tpfS$ of $\tfS$.

The \emph{height} of a point $(z,t)\in \widetilde{\fS}$ is $t$. A vector at $(z,t)$ is called \emph{vertical} if it is a positive vector of $\{z\} \times (-1,1)$. A \emph{framing} of a 1-dimensional submanifold $\al$ of $\fS \times (-1,1)$ is a continuous choice of a vector transverse to $\al$ at each point of $\al$. If $\al$ is equipped with a framing we say that $\al$ is \emph{framed}.

By a \emph{$\tpfS$-tangle} $\al$ in $\tfS= \fS\times (-1,1)$ we mean a framed 1-dimensional compact submanifold of $M$ such that
\begin{itemize}
\item $\partial \al \subset \tpfS$ and the framing at each point in $\partial \al$ is vertical, and
\item the boundary points of $\al$ in a boundary wall have distinct heights.
\end{itemize}
For any $\tpfS$-tangle $\al$, there is the partial order called the \emph{height order} on $\partial \al$, where $x>y$ if and only if $x,y$ are in the same boundary wall and the height of $x$ is greater than that of $y$. If $x>y$ and there is no $z$ such that $x>z>y$, then we say $x$ and $y$ are \emph{consecutive}.

Two $\tpfS$-tangles are \emph{isotopic} if they are isotopic in the class of $\tpfS$-tangles. In particular, $\tpfS$-isotopies do not change the height order. The empty set, by convention, is a $\tpfS$-tangle which is isotopic only to itself.

As usual, $\tpfS$-tangles are depicted by their diagrams on $\fS$ as follows. Every $\tpfS$-tangle is isotopic to one with vertical framing. Suppose a vertically framed $\tpfS$-tangle $\al$ is in general position with respect to the projection $\pi: \fS \times (-1,1) \to \fS$. The restriction $\pi |_{\al}:\al \to \fS$ is an immersion with transverse double points as the only possible singularities, and there are no double points on the boundary of $\fS$. Then $D=\pi(\al)$, together with
\begin{itemize}
\item the over/underpassing information at every double point, and
\item the linear order on $\pi(\al) \cap b$ for each boundary edge $b$ induced from the height order
\end{itemize}
is called a \emph{$\pfS$-tangle diagram}, or simply a tangle diagram on $\fS$. \emph{Isotopies} of $\pfS$-tangle diagrams are ambient isotopies in $\fS$.

Clearly the $\pfS$-tangle diagram of a $\tpfS$-tangle $\al$ determines the isotopy class of $\al$. When there is no confusion, we identify a $\pfS$-tangle diagram with its isotopy class of $\tpfS$-tangles.

Let $\ori$ be an orientation of $\pfS$, which may or may not be the orientation inherited from $\fS$. A $\pfS$-tangle diagram $D$ is \emph{$\ori$-ordered} if for each boundary edge $b$, the order of $\partial D$ on $b$ is increasing when one goes along $b$ in the direction of $\ori$. It is clear that every $\tpfS$-tangle can be presented by an $\ori$-ordered $\partial\surface$-tangle diagram after an isotopy. If $\ori$ is the orientation coming from $\fS$, the $\ori$-order is called the \emph{positive order}.

\subsection{Stated skein algebra}\label{sec.sk}

A \emph{state} on a finite set $X$ is a map $s: X \to \{\pm\}$. 
If $\al$ is a $\tpfS$-tangle or a $\pfS$-tangle diagram, then $\al$ is \emph{stated} if it is equipped with a state on the set $\partial \al$ of boundary points.

\newcommand{\relemp}{
\begin{tikzpicture}[scale=0.8,baseline=0.3cm]
\fill[gray!20!white] (0,0)rectangle(1,1);
\draw (1,0)--(1,1);
\end{tikzpicture}\,
}
\newcommand{\relconn}{
\begin{tikzpicture}[scale=0.8,baseline=0.3cm]
\fill[gray!20!white] (0,0)rectangle(1,1);
\draw (1,0)--(1,1);
\draw [very thick] (0,0.67)..controls(0.8,0.67) and (0.8,0.33)..(0,0.33);
\end{tikzpicture}\,
}
\newcommand{\relup}[2]{
\begin{tikzpicture}[scale=0.8,baseline=0.3cm]
\fill[gray!20!white] (0,0)rectangle(1,1);
\draw[-stealth] (1,0)--(1,1);
\draw[very thick] (0,0.67)--(1,0.67) (0,0.33)--(1,0.33);
\draw[inner sep=1pt,right] (1,0.67)node{\tiny #1} (1,0.33)node{\tiny #2};
\end{tikzpicture}
}
\newcommand{\reldown}[2]{
\begin{tikzpicture}[scale=0.8,baseline=0.3cm]
\fill[gray!20!white] (0,0)rectangle(1,1);
\draw[-stealth] (1,1)--(1,0);
\draw[very thick] (0,0.67)--(1,0.67) (0,0.33)--(1,0.33);
\draw[inner sep=1pt,right] (1,0.67)node{\tiny #1} (1,0.33)node{\tiny #2};
\end{tikzpicture}
}
\newcommand{\relarc}[2]{
\begin{tikzpicture}[scale=0.8,baseline=0.3cm]
\fill[gray!20!white] (0,0)rectangle(1,1);
\draw[-stealth] (1,0)--(1,1);
\draw[very thick] (1,0.67)..controls(0.2,0.67) and (0.2,0.33)..(1,0.33);
\draw[inner sep=1pt,right] (1,0.67)node{\tiny #1} (1,0.33)node{\tiny #2};
\end{tikzpicture}
}

Let $\cS(\fS)$ be the $\cR$-module freely spanned by isotopy classes of stated $\tpfS$-tangles modulo the \emph{defining relations}, which are the skein relation \eqref{eq.skein}, the trivial loop relation \eqref{eq.loop}, the trivial arc relations \eqref{eq.arcs}, and the state exchange relation \eqref{eq.order}:
\begin{align}
\label{eq.skein}
\begin{tikzpicture}[scale=0.8,baseline=0.3cm]
\fill[gray!20!white] (-0.1,0)rectangle(1.1,1);
\begin{knot}[clip width=8,background color=gray!20!white]
\strand[very thick] (1,1)--(0,0);
\strand[very thick] (0,1)--(1,0);
\end{knot}
\end{tikzpicture}
&=q
\begin{tikzpicture}[scale=0.8,baseline=0.3cm]
\fill[gray!20!white] (-0.1,0)rectangle(1.1,1);
\draw[very thick] (0,0)..controls (0.5,0.5)..(0,1);
\draw[very thick] (1,0)..controls (0.5,0.5)..(1,1);
\end{tikzpicture}
+q^{-1}
\begin{tikzpicture}[scale=0.8,baseline=0.3cm]
\fill[gray!20!white] (-0.1,0)rectangle(1.1,1);
\draw[very thick] (0,0)..controls (0.5,0.5)..(1,0);
\draw[very thick] (0,1)..controls (0.5,0.5)..(1,1);
\end{tikzpicture}\, ,\\
\label{eq.loop}
\begin{tikzpicture}[scale=0.8,baseline=0.3cm]
\fill[gray!20!white] (0,0)rectangle(1,1);
\draw[very thick] (0.5,0.5)circle(0.3);
\end{tikzpicture}
&=(-q^2 -q^{-2})
\begin{tikzpicture}[scale=0.8,baseline=0.3cm]
\fill[gray!20!white] (0,0)rectangle(1,1);
\end{tikzpicture}\, ,\\
\label{eq.arcs}
\relarc{$+$}{$-$}&=q^{-1/2}\relemp,\quad
\relarc{$+$}{$+$}=0, \quad \relarc{$-$}{$-$}= 0, \\
\label{eq.order}
\relup{$-$}{$+$}&=q^2\relup{$+$}{$-$}+q^{-1/2}\relconn.
\end{align}
Here each shaded part is a part of $\fS$, with a stated $\pfS$-tangle diagram on it indicated by the thick lines. Each thin line is part of a boundary edge, and the height order on that part is indicated by the arrow, and the points on that part are consecutive. The order of other endpoints away from the picture can be arbitrary and are not determined by the arrows of the pictures. 

For two $\tpfS$-tangles $\al_1$ and $\al_2$, the product $\al_1 \al_2$ is defined as the result of stacking $\al_1$ above $\al_2$. The product makes $\cS(\fS)$ an $\cR$-algebra, which is non-commutative is general.

Given a proper embedding of punctured bordered surfaces $\iota:\surface\to\surface'$ where each boundary edge of $\fS'$ contains the image of at most one boundary edge of $\fS$,
there is an induced map on the stated skein algebras
\[\iota_\ast:\stateS(\surface)\to\stateS(\surface')\]
defined on the stated $\partial\surface$-tangle diagrams in the obvious way. It is clear that $\iota_\ast$ respects
the product structure. Thus $\iota_\ast$ is a well defined $\mathcal{R}$-algebra homomorphism.

\begin{remark} \label{rem.hS}
Relations \eqref{eq.arcs} appeared in \cite{BW1}. Relation \eqref{eq.order} appeared in \cite{Le:TDEC} where the stated skein algebra was introduced.
\end{remark}

\def\tSS{\cS^{\mathrm{BW}}(\fS)}
\def\BW{{\mathrm{BW}}}

\subsection{Related constructions}

If we use only the relations \eqref{eq.skein}~and~\eqref{eq.loop} in the definition of $\SS$, we get a coarser version $\tSS$, which was defined by Bonahon and Wong \cite{BW1}.

The subalgebra $\mullerS(\surface)$ spanned by $\partial\widetilde{\surface}$-tangles whose states are all $+$ is naturally isomorphic to the skein algebra defined by Muller \cite{Muller}. For more details, see \cite{Le:TDEC,LP}.

The ordinary skein algebra is defined using links (tangles with only closed components) and the relations \eqref{eq.skein}~and~\eqref{eq.loop}. Since the interior $\ofS$ of $\fS$ does not have boundary, arc components are not allowed in the thickening of $\mathring{\surface}$, and the relations \eqref{eq.arcs}~and\eqref{eq.order} are vacuous. Thus $\cS(\ofS)$ is identified with the ordinary skein algebra of $\surface$. The inclusion $\mathring{\surface}\hookrightarrow\surface$ induces a natural algebra embedding $\cS(\mathring{\surface}) \embed \SS$, see \cite{Le:TDEC}.

\subsection{Height exchange relations, reflection}

\newcommand{\relcross}[2]{
\begin{tikzpicture}[scale=0.8,baseline=0.3cm]
\fill[gray!20!white] (0,0)rectangle(1,1);
\draw[-stealth] (1,0)--(1,1);
\begin{knot}[clip width=4,background color=gray!20!white]
\strand[very thick] (0,0.3)--(1,0.67);
\strand[very thick] (0,0.7)--(1,0.33);
\end{knot}
\draw[inner sep=1pt,right] (1,0.67)node{\tiny #1} (1,0.33)node{\tiny #2};
\end{tikzpicture}
}

\begin{lemma}[Height exchange relations, Lemma~2.4 of \cite{Le:TDEC}]
\label{lemma-height}
For $\nu\in \{\pm \}$ one has
\begin{align}
\label{eq.reor1}
\relcross{$\nu$}{$+$} = q^{-\nu}\relup{$+$}{$\nu$}, \quad
\relcross{$-$}{$\nu$} &= q^\nu\relup{$\nu$}{$-$}, \\
\label{eq.reor2}
q^{3/2}\reldown{$-$}{$+$}-q^{-3/2}\relup{$-$}{$+$} &= (q^2-q^{-2})\relconn.
\end{align}
Here we have identified $\pm$ with $\pm 1$ when we write $q^\nu$.
\end{lemma}


\begin{proposition}[Reflection anti-involution, Proposition 2.7 in \cite{Le:TDEC}] \label{lemma-r.reflection}
When $\cR=\BZ[q^{\pm 1/2}]$, there exists a unique $\BZ$-linear map $\omega: \cS(\fS) \to \cS(\fS)$ such that
\begin{itemize}
\item $\omega(q^{1/2})= q^{-1/2}$,
\item $\omega$ is an anti-automorphism, and $\omega^2=\id$,
\item if $\al$ is a stated $\pfS$-tangle diagram then $\omega(\al)$ is the result of switching all the crossings of $\al$ and reversing the height order on each boundary edge.
\end{itemize}
\end{proposition}
We call $\omega$ the \emph{reflection anti-involution}. A map $f:\stateS(\surface)\to\stateS(\surface')$ is \emph{reflection invariant} if $f\circ\omega=\omega\circ f$. For example, the map $\iota_\ast$ induced by a surface embedding $\iota$ is reflection invariant.

Similarly, for an antisymmetric integer matrix $Q$, a map $f: \CS\to \bT(Q)$ is \emph{reflection invariant} if $f\circ\omega=\omega\circ f$.

\subsection{Basis} \label{sec.basis0}

A $\pfS$-tangle diagram $\alpha$ is \emph{simple} if it has neither double point nor trivial component. Here a closed component of $\alpha$ is \emph{trivial} if it bounds a disk in $\fS$, and an arc component of $\al$ is \emph{trivial} if it can be homotoped relative to its boundary to a subset of a boundary edge. By convention, the empty set is considered as a simple stated $\partial\widetilde{\surface}$-tangle diagram.

Define an order on $\{\pm\}$ so that the sign $-$ is less than the sign $+$. Recall the set of boundary points $\partial\alpha$ has a partial order (the height order). A state $s: \partial\alpha\to \{\pm\}$ is \emph{increasing} if $s$ is an increasing function, i.e., $s(x) \le s(y)$ whenever $x \le y$.

Let $B(\fS)$ be the set of of all isotopy classes of increasingly stated, positively ordered simple $\partial\surface$-tangle diagrams.

\begin{theorem}[Theorems 2.11 and Proposition 4.4 in \cite{Le:TDEC}] \label{thm.basis}
Suppose $\fS$ is a punctured bordered surface. The set $B(\fS)$ is a free $\cR$-basis of $\cS(\fS)$.
The algebra $\CS$ is a domain. In other words, if $xy =0$ where $x,y\in \CS$, then $x=0$ or $y=0$.
\end{theorem}

\def\bbn{\bar{\bn}}

\def\tphi{\tilde \phi}
\def\fB{{\mathfrak{B}}}
\def\SB{\cS(\fB)}
\def\SrB{\cS^\rd(\fB)}
\def\YB{\cY(\fB)}
\def\he{{\hat e}}

\subsection{Splitting homomorphism}\label{sec.splitting}
Suppose $c$ is an ideal arc in the interior of a puncture bordered surface $\fS$. The \emph{splitting} of $\fS$ along $c$ is a punctured bordered surface $\fS'$ with two boundary edges $a,b$ such that gluing $a$ and $b$ in $\surface'$ gives $\surface$ and $a$ and $b$ both project to $c$.

A $\tpfS$-tangle $\al\subset \tfS$ is \emph{vertically transverse to $c$} if
\begin{itemize}
\item $\al$ is transverse to $c \times (-1,1)$, and
\item the points in $\partial_c\al \eqdef \al \cap (c \times (-1,1))$ have distinct heights and vertical framing.
\end{itemize}
Suppose $\al$ is a stated $\tpfS$-tangle vertically transverse to $c$. By splitting $\al$ along $c\times (-1,1)$ we get a $\tpfS'$-tangle $\tal$ which is stated at every boundary point except for points in $p^{-1}(\partial_c\al)$, where $p:\fS' \times (-1,1) \to \tfS=\fS\times (-1,1)$ is the natural projection. For every $s: \partial_c\al \to \{\pm \}$ let $(\tal, s)$ be the $\tpfS'$-tangle $\tal$ where the state of a point $x$ in $p^{-1}(\partial_c\al)$ is $s(p(x))$.

\begin{theorem}[Splitting Theorem, Theorem~3.1 in \cite{Le:TDEC}] \label{thm.1a} Suppose $c$ is an ideal arc in the interior of a punctured bordered surface $\fS$, and $\fS'$ is the splitting of $\fS$ along $c$.
\begin{enumerate}
\item There is a unique $\cR$-algebra homomorphism $\theta_c: \cS(\fS) \to \cS(\fS')$, called the splitting homomorphism along $c$, such that if $\al$ is a stated $\tpfS$-tangle vertically transverse to $c$, then
\begin{equation}\label{eq.split}
\theta_c(\al)=\sum_{s: \partial_c \al \to \{\pm \}} (\tal,s).
\end{equation}
\item In addition, $\theta_c$ is injective and reflection invariant.
\item If $c_1$ and $c_2$ are two disjoint ideal arcs in the interior of $\fS$, then
\[\theta_{c_1} \circ \theta_{c_2} = \theta_{c_2} \circ \theta_{c_1}.\]
\end{enumerate}
\end{theorem}

\begin{remark}
If $\fS_1$ and $\fS_2$ are two punctured bordered surfaces, then there is a natural algebra isomorphism
\[\cS (\fS_1) \ot_\cR \cS (\fS_2)\cong \cS (\fS_1 \sqcup \fS_2)\]
where $\alpha_1\otimes\alpha_2$ is identified with $\alpha_1\sqcup\alpha_2$ for $\partial\widetilde{\surface}_i$-tangles $\alpha_i$, $i=1,2$. Thus if the splitting of $\fS$ along $c$ is the disconnected surface $\surface_1\sqcup\surface_2$, we can write the splitting homomorphism as \[\theta_c:\stateS(\surface)\to\stateS(\surface_1)\otimes\stateS(\surface_2).\]
This statement easily generalizes to surfaces with more than two components.
\end{remark}

\def\tbfS{\widetilde{\bfS}}
\def\tPd{\widetilde{\Pd}}

\subsection{Bigon and coaction on surface} \label{sec.bigon}

\begin{figure}[h]
\centering
\begin{tikzpicture}
\draw[fill=gray!20!white] (0,-1)
	arc[start angle=-60,end angle=60,radius=1.15]
	arc[start angle=120,end angle=240,radius=1.15];
\draw[fill=white] (0,-1)circle(2pt) (0,1)circle(2pt);
\draw[very thick] (-0.58,0)node[left]{$\mu$} to node[above]{$a$}
	(0.58,0)node[right]{$\nu$};
\draw (0.64,0.7)node{$e_R$} (-0.64,0.7)node{$e_L$};
\end{tikzpicture}
\caption{Bigon and arc $a_{\mu\nu}$.}\label{fig:bigon}
\end{figure}
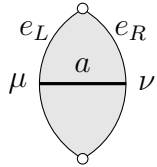

Let $\cB$ be an ideal bigon, with two boundary edges $e_L$ and $e_R$. For $\mu, \nu\in \{\pm \}$ let $a_{\mu\nu}$ be the stated arc depicted in Figure~\ref{fig:bigon}. The algebra $\SB$ is generated by $a_{\mu\nu}$ with $\mu, \nu\in \{\pm \}$. There are geometrically defined coproduct, counit, and antipode for $\SB$ which turns $\SB$ into a Hopf $\cR$-algebra isomorphic to the well-known quantum coordinate ring $\OSL$ of the Lie group $SL_2$, see \cite{CL,KQ}.

The counit $\ve : \SB \to \cR$ is the $\cR$-algebra homomorphism defined on generators by
\begin{equation}\label{eqn-counit}
\ve(a_{\mu\nu}) = \delta_{\mu,\nu}.
\end{equation}

If $\al$ is a stated $\partial \tilde \fB$-tangle, then $\ve(\al)$ is the matrix element of the Reshetikhin-Turaev operator invariant of the tangle $\al$, see \cite{CL}, and is also equal to $\operatorname{Tr}_\mathfrak{B}(\al)$ of \cite[Proposition 13]{BW1}. There is a conservation of charge property, which says for any $\partial \fB$-tangle diagram $\al$ with state $s: \partial \al \to \{\pm \} \equiv \{\pm 1 \}$,
\begin{equation}
\ve(\al) = 0\text{ if } \sum_{x\in \partial\alpha \cap e_L} s(x) \neq \sum_{x\in \partial\alpha \cap e_R} s(x).
\label{eq.charge}
\end{equation}

\def\ST{\cS(\fT)}
\def\YT{\cY(\fT)}

\begin{lemma}\label{r.plus}
Suppose $\al$ is a stated simple $\partial \fB$-tangle diagram whose states are positive. Then $\ve(\al)=q^{k/2}$ for some $k\in\mathbb{Z}$. In particular, it is nonzero.
\end{lemma}

\begin{proof}
Use a diagram $D$ of $\al$ where the height orders on $e_L$ and $e_R$ are given from bottom to top. In general $D$ has crossing. Because the states are all positive, there is only one way to resolve all the crossings of $D$ so that the resulting diagram is not 0. Thus $D$ is $q$-proportional to the diagram with horizontal parallel arcs, whose $\ve$ value is $1$ by \eqref{eqn-counit}.
\end{proof}

Suppose $e$ is a boundary edge of a punctured bordered surface $\fS$. By splitting $\fS$ along an ideal arc $e'$ isotopic to $e$ and lying in the interior of $\fS$, we get a bigon and a surface diffeomorphic to $\fS$. The splitting homomorphism gives an algebra map
\begin{equation}
\CS \to \CS \ot \SB
\end{equation}
which provides $\CS$ with a right $\stateS(\cB)$-comodule structure. See \cite{CL} for properties of this comodule. In particular, when $\fS=\cB$, the above comodule map is the coproduct.

\section{Embedding into quantum torus, length coordinate version}
\label{sec.length}

Throughout this section $\fS$ is a connetced punctured bordered surface with non-empty boundary. In particular, the set $\Pd$ of boundary punctures is non-empty. For simplicity we also assume that $\fS$ is not a monogon or a bigon.

We will show that there is an embedding of $\cS(\fS)$ into a quantum torus
\[\varphi: \CS \embed \bT(\bar{P})\]
such that the image of $\varphi$ is sandwiched $\bT_+(\bar{P})$ and $\bT(\bar{P})$.
This fact makes it easy to study the representations of $\CS$. The matrix $\bar{P}$ depends on a \emph{quasitriangulation} of $\fS$.

When $\fS$ does not have interior punctures and $q=1$, and $\al$ is a simple closed curve, the image $\varphi(\al)$ expresses the lambda length of $\al$ as a Laurent polynomial in the lambda length coordinates of the decorated Teichm\"uller space.

\subsection{Disjoint ideal arcs}\label{sec.P.arc}

\newcommand{\vertexdef}[2]{
\begin{tikzpicture}[scale=0.8,baseline=0.28cm]
\fill[gray!20!white] (-0.2,0) rectangle (1.3,1);
\draw (0,1)--(0.6,0)--(1,1) (-0.2,0)--(1.3,0);
\draw (0.1,0.5)node{\vphantom{$b$}#1} (0.98,0.5)node{\vphantom{$b$}#2};
\draw[fill=white] (0.6,0)circle(2pt);
\end{tikzpicture}
}

An ideal arc of $\fS$ is \emph{boundary ending} if its endpoints are boundary punctures. For disjoint boundary ending ideal arcs $a$ and $b$ (which can be isotopic), define the integer $P(a,b)$ by
\[P(a,b)= \#\left( \vertexdef{$b$}{$a$} \right) - \#\left( \vertexdef{$a$}{$b$} \right).\]
To be precise, by removing a point in $a$ we get two \emph{half-edges of $a$}, each is incident with exactly one boundary puncture. For a boundary puncture $v$ and two disjoint half-edges $a', b'$ let $P_v(a',b')=0$ if one of $a', b'$ is not incident with $v$; otherwise let
\[P_v(a',b')=\begin{cases}
1,&b'\text{ is counterclockwise to }a',\\
-1,&b'\text{ is clockwise to }a'.
\end{cases}\]
Let $P(a,b)=\sum P_v(a',b')$ where the sum is over half-edges $a'$ of $a$ and half-edges $b'$ of $b$ and all boundary punctures $v$. It is evident that $P(a,b)$ only depends on the isotopy classes of $a$ and $b$.

\subsection{Quasitriangulation and its vertex matrix}
An \emph{ideal multiarc} in $\fS$ is a finite collection of disjoint ideal arcs, and it is \emph{boundary ending} if each component is boundary ending. In this paper, a \emph{quasitriangulation} of $\fS$ is a maximal boundary ending ideal multiarc $\cE$ whose components are non-trivial and no two of them are isotopic. Any boundary edge of $\fS$ is isotopic to an element of $\cE$. Since $\fS$ is not a monogon nor a bigon, every boundary edge of $\fS$ is isotopic to an element of $\cE$, and two different boundary edges are not isotopic. Thus, after an isotopy, we can assume that $\cE$ contains the set $\Ed$ of all boundary edges. Let $\mathring{\quasi}=\quasi\setminus\quasi_\partial$ be the set of \emph{interior ideal arcs} of $\cE$. If we split $\fS$ along all $e\in \mathring{\quasi}$, we get a collection of ideal triangles and once-punctured monogons.

Let $P:\cE\times \cE\to \BZ$ be the antisymmetric function whose values $P(a,b)$ are defined in Section~\ref{sec.P.arc}.
Let $\hEd= \{ \hat e \mid e\in \Ed\}$ be another copy of the set $\Ed$ of boundary edges, and let $\bcE= \cE \sqcup \hEd$. Define an antisymmetric function $\bP: \bcE \times \bcE \to \BZ$, which is an extension of $P$, by
\begin{alignat}{2}
\bP(a,\hat b) &= -\#\left( \vertexdef{$b$}{$a$} \right) -\#\left( \vertexdef{$a$}{$b$} \right),& &\text{if }a\in \cE, b\in \Ed, \label{eq.P2}\\
\bP(\hat a, \hat b) &= -P(a,b), &\quad&\text{if }a,b \in \Ed. \notag
\end{alignat}
The right hand side of \eqref{eq.P2} counts the number of times when a half-edge of $a$ and a half-edge of $b$ meet at a boundary puncture (with a minus sign).

Recall that $\Po$ is the set of interior punctures. Let $\cR[\Po]$ be the polynomial algebra in variables which are the interior punctures. Define the quantum torus associated to the antisymmetric form
$\bP$ with ground ring $\cR[\Po]$:
\begin{align*}
\bX(\fS;\cE) &= \cR[\Po]\la x_a^{\pm 1}, a \in \bcE\ra /( x_a x_b = q^{\bar{P}(a,b)} x_b x_a).
\end{align*}
The set $\{ x^\bk\mid \bk \in \BZ^{\bcE}\}$ is a free $\cR[\Po]$-basis of $\bX(\fS;\cE)$. The $\cR[\Po]$-submodule spanned by $\{ x^\bk\mid \bk \in \BN^{\bcE}\}$ is denoted by $\bX_+(\fS;\cE)$.

\subsection{Quantum trace, length coordinates version}

\begin{theorem}\label{thm.embed1}
Suppose $\fS$ is a connected punctured bordered surface with non-empty boundary, and $\fS$ is not a monogon or a bigon. Let $\cE$ be a quasitriangulation of $\fS$.

There is an $\cR$-algebra embedding $\varphi_\cE:\cS(\fS) \embed \bX(\fS;\cE)$ such that
\begin{equation}
\bX_+(\fS;\cE)\subset\varphi_\quasi(\cS(\fS)) \subset \bX(\fS;\cE).
\end{equation}
Consequently $\cS(\fS)$ is an Ore domain, and
$\vpE$ induces an isomorphism of the division rings $\tilde{\varphi}_\quasi:\Fr(\SS) \xrightarrow{\cong} \Fr( \bX(\fS;\cE) )$.

When $\cR=\BZ[q^{\pm 1/2}]$, the map $\vpE$ is reflection invariant.
\end{theorem}

The fact that $\varphi(\cS(\fS))$ is sandwiched between the quantum space $\bX_+(\fS;\cE)$ and the quantum torus $\bX(\fS;\cE)$ has many applications. For example, it follows that the center of $\CS$ is the restriction of the center of $\extX(\surface;\quasi)$. Any representation of $\CS$ in which the actions of $x_a, a\in \cE$ are invertible extends to a representation of the quantum torus whose representation theory is known.

Another immediate consequence is the easy construction of the coordinate change isomorphism. For another quasitriangulation $\cE'$, define the \emph{coordinate change isomorphism}
\[\Psi_{\cE' \cE}: \Fr(\bX(\fS;\cE)) \to \Fr(\bX(\fS;\cE')), \quad
\Psi_{\cE' \cE}= \tilde{\varphi}_{\cE'} \circ (\tilde{\varphi}_\cE)^{-1}.\]
By construction, $\Psi_{\quasi\quasi'}$ is reflection invariant, that is, $\Psi_{\quasi\quasi'}$ commutes with the extension $\tilde{\omega}$ of the reflection anti-involution to the ring of fractions.

\begin{corollary}
The coordinate change isomorphism is functorial in the sense that for quasitriangulations $\cE, \cE', \cE''$, one has
\[\Psi_{\cE \cE}= \id,\quad \Psi_{\cE'' \cE} = \Psi_{\cE'' \cE'}\circ \Psi_{\cE' \cE}.\]
In addition, $\varphi_{\cE'}= \Psi_{\cE' \cE} \circ \varphi_{\cE}$.
\end{corollary}

It should be noted that the coordinate change isomorphism for the (quantum) coordinate embedding of Section~\ref{sec-shear} is much more difficult to construct. In \cite{Liu,Hiatt,BW1}, one first constructed the coordinate change isomorphism (for the shear quantum trace map) for flips of triangulations, then showed that the isomorphism does not depend on the sequence of flips representing the change of triangulations.

\begin{remark}
When $\fS$ does not have interior puncture, M\"uller \cite{Muller} first defined the matrix $P$ and showed that there is an embedding of the subalgebra $\cS^+(\fS)$ into the smaller quantum torus $\sX(\fS;\cE) \subset \bX(\fS;\cE)$ generated by $x_a^{\pm 1}, a\in \cE$. He showed that when $q=1$, the embedding expresses the lambda length of a simple closed curve as a Laurent polynomial in the lambda lengths of the edges of the quasitriangulation. The embedding result is extended to the case when $\fS$ has interior punctures in \cite{LP}. Here we extend the result to the whole algebra $\cS(\fS)$, using the bigger the quantum torus $\bX(\fS;\cE)$.
\end{remark}

\def\he{{\hat e}}
\def\RV{{\cR[\Po]}}
\def\Bo{\mathbullet B}
\subsection{Ideal arcs as elements of $\SS$}\label{sec.bad}


The polynomial ring $\RV$ can be considered as a subalgebra of the center of $\SS$, where we identify $v\in \Po$ with a small loop $X_v$ in $\fS$ surrounding $v$. Such a loop will be called a \emph{peripheral loop}. From Theorem~\ref{thm.basis}, the set $\Bo$ consisting of positively ordered, increasingly stated, simple $\pfS$-tangle diagrams without peripheral loops is a free $\RV$-basis of $\SS$.


\begin{figure}[h]
\centering
\begin{tikzpicture}
\begin{scope}
\clip (1,0)arc[x radius=1.25,y radius=1,start angle=0,end angle=180];
\fill[gray!20!white] (-1.5,0) rectangle (1,1);
\draw[very thick] (0,0) edge (0.5,1) edge (-0.5,1) edge (-1.5,1);
\end{scope}
\draw (-1.5,0)--(1,0) (-0.5,0)node[below]{$\vphantom{D}\alpha$};
\draw[fill=white] (0,0)circle(2pt);

\begin{scope}[xshift=3.5cm]
\begin{scope}
\clip (1,0)arc[x radius=1.25,y radius=1,start angle=0,end angle=180];
\fill[gray!20!white] (-1.5,0) rectangle (1,1);
\draw[very thick] (0.5,1)--(-0.2,0) (-0.5,1)--(-0.4,0) (-1.5,1)--(-0.6,0);
\end{scope}
\path[tips,-{Stealth[length=0.2cm]}] (-1.5,0)--(-1,0);
\draw (-1.5,0)--(1,0) (-0.5,0)node[below]{$D(\alpha)$};
\draw[fill=white] (0,0)circle(2pt);
\end{scope}
\end{tikzpicture}
\caption{Moving left: From ideal multiarc $\al$ to $D(\al)$.}\label{fig:movingleft}
\end{figure}
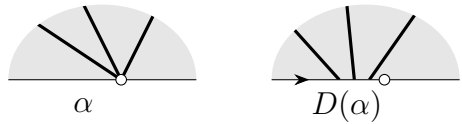

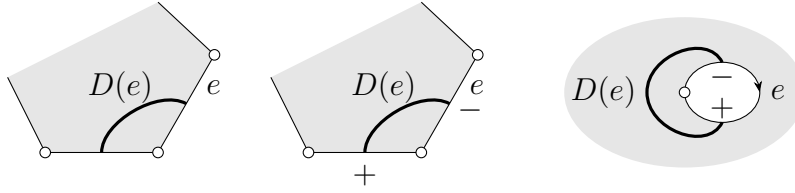
\begin{figure}[h]
\centering
\begin{tikzpicture}
\draw[fill=gray!20!white] (-2,1)--(-1.5,0)--(0,0)--(60:1.5)--(0,2);
\draw[very thick] (60:0.75)..controls (0,0.85) and (-0.75,0.4)..(-0.75,0);
\draw (120:1)node{$D(e)$} (60:1)node[right]{$e$};
\draw[fill=white] (0,0)circle(2pt) (-1.5,0)circle(2pt) (60:1.5)circle(2pt);
\begin{scope}[xshift=3.5cm]
\draw[fill=gray!20!white] (-2,1)--(-1.5,0)--(0,0)--(60:1.5)--(0,2);
\draw[very thick] (60:0.75)..controls (0,0.85) and (-0.75,0.4)..
	(-0.75,0)node[below]{$+$} (60:0.65)node[right]{$-$};
\draw (120:1)node{$D(e)$} (60:1)node[right]{$e$};
\draw[fill=white] (0,0)circle(2pt) (-1.5,0)circle(2pt) (60:1.5)circle(2pt);
\end{scope}
\begin{scope}[xshift=7cm,yshift=0.8cm]
\fill[gray!20!white] (0,0)circle[x radius=1.6,y radius=1];
\draw[fill=white] (0.5,0)circle[x radius=0.5,y radius=0.4];
\path[tips,-{Stealth[length=0.15cm]}] (0.82,1)--(1,0);
\draw[very thick,inner sep=1pt] (0.5,0.4)node[below]{$-$}
	..controls (0.3,0.8) and (-0.5,0.5)..(-0.5,0)
	..controls (-0.5,-0.5) and (0.3,-0.8)..(0.5,-0.4)node[above]{$+$};
\draw[fill=white] (0,0)circle(2pt) (1,0)node[right]{$e$} (-0.5,0)node[left]{$D(e)$};
\end{scope}
\end{tikzpicture}
\caption{Left: Boundary edge $e$ and corner arc $D(e)$. Middle: Bad arc $D(e)(+,-)$. Right: Bad arc $D(e)(+,-)$ when two endpoints of $e$ coincide.}\label{fig:barc1}
\end{figure}

For each ideal multiarc $\al$, let $D(\al)$ be the simple $\pfS$-tangle diagram obtained from $\al$ by slightly moving all the strands of $\al$ coming to each boundary puncture $p$ to strands ending on the boundary edge lying to the left of $p$, see Figure~\ref{fig:movingleft}, and imposing positive order on each boundary edge. Note that $D(\al)$ is well-defined up to isotopy. When $e\in\Ed$ is a boundary edge $D(e)$ is called a \emph{corner arc}, see Figure~\ref{fig:barc1}.

For each $e\in \bcE$ define $X_e\in \cS(\fS)$ as follows. For $e\in \cE$ let
\begin{align*}
X_e= \begin{cases}
D(e)(+,+), & \text{if endpoints of $e$ are disctinct},\\
q^{-1/2} D(e)(+,+), & \text{if endpoints of $e$ are identical},
\end{cases}
\end{align*}
where $D(e)(+,+)$ is $D(e)$ equipped with positive states at both endpoints. The prefactor $q^{-1/2}$ is introduced so that $X_e$ is reflection invariant (when $\cR=\BZ[q^{\pm 1/2}]$), which follows from the height exchange relation~\eqref{eq.reor1}.

For $e\in \Ed$ define $X_\he\in \SS$, called the \emph{bad arc} corresponding to $e$ following \cite{CL}, by
\begin{align*}
X_\he= \begin{cases}
D(e)(+,-), & \text{if endpoints of $e$ are disctinct},\\
q^{1/2} D(e)(+,-), & \text{if endpoints of $e$ are identical},
\end{cases}
\end{align*}
where $D(e)(+,-)$ is $D(e)$ with states assigned as in Figure~\ref{fig:barc1}. Again $X_\he$ is reflection invariant.

\begin{lemma} \label{r.badarc}
Let $\alpha$ be a stated $\pfS$-tangle diagram.
\begin{enumerate}
\item In $\cS(\fS)$, $\al$ is $q$-commuting with any bad arc
$X_\he$. Consequently $X_\he$ is a normal element of $\cS(\fS)$ in the sense that the left ideal generated by $X_\he$ is the same as the right ideal generated by $X_\he$.

\item Suppose $\al\in \Bo$. Let $\al'$ be the union of all the bad arcs in $\al$ and $\al''$ is the union of the remaining components. Then $\al' \al''$ is $q$-proportional to $\al$.
\end{enumerate}
\end{lemma}

\begin{proof}
For part (a), first assume that the two endpoints of $e$ are distinct. Let $a$ be the boundary edge lying to the left of $e$. The diagrams of the products $X_\he \al$ and $\al X_\he$, together with the diagram $Y$ which is a disjoint union of $\alpha$ and $X_{\hat{e}}$, are given in Figure~\ref{fig:badarc}. All boundary orders are positive.

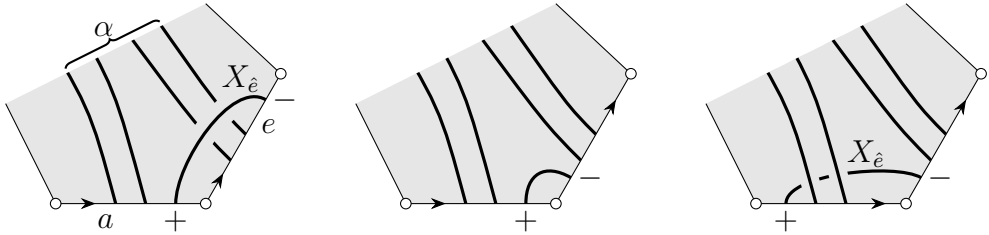
\begin{figure}[h]
\centering
\begin{tikzpicture}[scale=1.33]
\fill[gray!20!white] (-2,1)--(-1.5,0)--(0,0)--(60:1.5)--(0,2);
\begin{scope}
\clip (-2,1)--(-1.5,0)--(0,0)--(60:1.5)--(0,2);
\draw[very thick] (-0.9,0)..controls ++(105:1)..(-1.8,2)
	(-0.6,0)..controls ++(105:1.2)..(-1.4,2)
	(60:0.5)..controls ++(135:0.5)..(-1,2)
	(60:0.8)..controls ++(135:0.4)..(-0.6,2);
\end{scope}
\draw[knot,very thick,background color=gray!20!white,knot gap=8,inner sep=2pt]
	(60:1.2)node[right]{$-$}..controls ++(150:0.4) and (-0.3,0.4)
	..(-0.3,0)node[below]{$+$};
\draw (-2,1)--(-1.5,0)--(0,0)--(60:1.5)--(0,2);
\path[tips,-{Stealth[length=0.2cm]}] (-1.5,0)--(-1.1,0);
\path[tips,-{Stealth[length=0.2cm]}] (0,0)--(60:0.35);
\draw[fill=white] (0,0)circle(1.5pt) (-1.5,0)circle(1.5pt) (60:1.5)circle(1.5pt);
\draw (60:1.5)++(-0.4,-0.05)node{$X_{\hat{e}}$} (60:0.9)node[right]{$e$} (-1,0)node[below]{$a$};
\draw[thick,decorate,decoration=brace] (-1.4,1.35)--(-0.46,1.82);
\draw (-0.95,1.6)++(-0.07,0.14)node{$\alpha$};

\begin{scope}[xshift=3.5cm]
\fill[gray!20!white] (-2,1)--(-1.5,0)--(0,0)--(60:1.5)--(0,2);
\begin{scope}
\clip (-2,1)--(-1.5,0)--(0,0)--(60:1.5)--(0,2);
\draw[very thick] (-0.9,0)..controls ++(105:1)..(-1.8,2)
	(-0.6,0)..controls ++(105:1.2)..(-1.4,2)
	(60:0.5)..controls ++(135:0.5)..(-1,2)
	(60:0.8)..controls ++(135:0.4)..(-0.6,2);
\end{scope}
\draw[very thick,inner sep=2pt] (60:0.3)node[right]{$-$}
	..controls ++(150:0.3) and (-0.3,0.3)
	..(-0.3,0)node[below]{$+$};
\draw (-2,1)--(-1.5,0)--(0,0)--(60:1.5)--(0,2);
\path[tips,-{Stealth[length=0.2cm]}] (-1.5,0)--(-1.1,0);
\path[tips,-{Stealth[length=0.2cm]}] (0,0)--(60:1.2);
\draw[fill=white] (0,0)circle(1.5pt) (-1.5,0)circle(1.5pt) (60:1.5)circle(1.5pt);
\end{scope}

\begin{scope}[xshift=7cm]
\fill[gray!20!white] (-2,1)--(-1.5,0)--(0,0)--(60:1.5)--(0,2);
\draw[very thick,inner sep=2pt] (60:0.3)node[right]{$-$}
	..controls ++(150:0.3) and (-1.2,0.3)
	..(-1.2,0)node[below]{$+$};
\begin{scope}
\clip (-2,1)--(-1.5,0)--(0,0)--(60:1.5)--(0,2);
\draw[knot,very thick,background color=gray!20!white,knot gap=6]
	(-0.9,0)..controls ++(105:1)..(-1.8,2)
	(-0.6,0)..controls ++(105:1.2)..(-1.4,2)
	(60:0.5)..controls ++(135:0.5)..(-1,2)
	(60:0.8)..controls ++(135:0.4)..(-0.6,2);
\end{scope}
\draw (-2,1)--(-1.5,0)--(0,0)--(60:1.5)--(0,2);
\draw (-0.4,0.5)node{$X_{\hat{e}}$};
\path[tips,-{Stealth[length=0.2cm]}] (-1.5,0)--(-0.2,0);
\path[tips,-{Stealth[length=0.2cm]}] (0,0)--(60:1.2);
\draw[fill=white] (0,0)circle(1.5pt) (-1.5,0)circle(1.5pt) (60:1.5)circle(1.5pt);
\end{scope}
\end{tikzpicture}
\caption{From left to right: Diagram of $X_\he \al$, Diagram $Y$, and Diagram of $\al X_\he$.}\label{fig:badarc}
\end{figure}

For a boundary edge $c$ let $\delta_c(\al)$ be the sum of all the states of boundary points of $\al$ on $c$, where we identify $\pm$ with $\pm 1$. The height exchange relation \eqref{eq.reor1} gives
\begin{align}
X_\he \al & = q^{\delta_e(\al)} Y \label{eq.bad1}\\
&= q^{\delta_e(\al) + \delta_a(\al) } \al X_\he.
\end{align}
When the two endpoints of $e$ are identical, then $a=e$ and a similar calculation shows
\begin{equation}
X_\he \al = q^{2\delta_e(\al)} \al X_\he.
\end{equation}

If $\alpha$ is positively ordered and increasingly stated, then the bad arc components looks like diagram $Y$. Thus (b) follows Identify \eqref{eq.bad1}.
\end{proof}

\def\Ein{\mathring \quasi}

Let $\Bo^\bad\subset \Bo$ be the set of stated $\pfS$-tangle diagrams in $\Bo$ whose components are bad arcs only, and $\Bo^\rd \subset \Bo$ be the set of stated $\pfS$-tangle diagrams in $\Bo$ containing no bad arcs.
Let $\cD(\fS)$ be the $\RV$-subalgebra of $\SS$ generated by the bad arcs. Then $\Bo^\bad$ is a free $\RV$-basis of $\cD(\fS)$. From Lemma~\ref{r.badarc}(b) we have the following corollary, which actually appeared implicitly in \cite[Section 7.2]{CL}.

\begin{corollary} \label{r.badbasis}
As a left $\cD(\fS)$-module $\SS$ is freely spanned by $\Bo^\rd$.
\end{corollary}

\subsection{Proof of Theorem~\ref{thm.embed1}}\label{sec.P.proof}
\begin{proof}

From the height exchange relation \eqref{eq.reor1} one gets that for any $a,b\in \bcE$,
\[X_a X_b = q^{ \bP(a,b) } X_b X_a.\]
This means there is a well-define $\RV$-algebra homomorphism
\[\iota:\extX_+(\surface;\quasi)\to\stateS(\surface), \quad \iota (x_e) = X_e. \]
which is reflection invariant (when $\cR= \BZ[q^{\pm 1/2}]$).

Now we show that $\iota$ is injective by showing that $\iota$ maps the $\RV$-basis $\{x^{\bk} \mid \bk\in \BN^{\bcE}\}$ of $\extX_+(\surface;\quasi)$ injectively into a subset of an $\RV$-basis of $\cS(\fS)$.
For $\bk \in \BN^{\bcE}$ let $\bk'\in \BN^{\cE}$ be defined by $\bk'(e)= \bk(e)$ if $e\in \Ein$ and $\bk'(e)= \bk(e)+ \bk(\hat e)$ for $e\in \Ed$. Let $\cE^{\bk'}$ be the ideal multiarc consisting of $\bk'(e)$ parallel copies of $e$ for every $e\in \cE$. Let $D(\cE^{\bk'};\bk)$ be the $\pfS$-tangle diagram $D(\cE^{\bk'})$ equipped with the increasing states such that on each boundary edge $e$ there are exactly $\bk(\hat{e})$ negative states. From the height exchange relation \eqref{eq.reor1}, we have that $\iota(x^\bk)$ is $q$-proportional to $D(\cE^{\bk'};\bk)$. Note that $D(\cE^{\bk'};\bk)\in \Bo$, and if $\bk_1\neq \bk_2$ then $D(\cE^{\bk_1'};\bk_1) \neq D(\cE^{\bk_2'};\bk_2)$. It follows that $\iota$ is injective.

Let us identify $\extX_+(\surface;\quasi)$ with a subset of $\cS(\fS)$ by the embedding $\iota$.

\begin{lemma}\label{lemma-Laurent}
For any $\alpha\in\stateS(\surface)$, there $\bk\in \BN^{\bcE}$ such that $x^{\bk} \al \in \extX_+(\surface;\quasi)$. Moreover one can choose $\bk$ such that $\bk(\he)=0$ for all $e\in \Ed$.
\end{lemma}

Since $\SS$ is a domain (see Theorem \ref{thm.basis}),
Lemma \ref{lemma-Laurent} and
Proposition~\ref{r.Ore} show that the embedding $\extX_+(\surface;\quasi) \embed \extX(\surface;\quasi)$ has a unique extension which is an $\RV$-algebra embedding
$\varphi_\cE: \cS(\fS) \embed \extX(\surface;\quasi)$, proving the theorem.

It remains to prove Lemma~\ref{lemma-Laurent}, which is a generalization \cite[Corollary 6.9]{Muller} and \cite[Lemma 6.5]{LP}, where it was showed that the lemma is true when $\al$ has only positive states.

We use induction on the number of negative states of $\al$, denoted $d(\al)$. If $d(\al)=0$, then $\al$ is positively stated, and the previous results apply. Now assume $d(\al)>0$.

Recall that $\cD(\fS)$, defined in Subsection~\ref{sec.bad} as the $\RV$-subalgebra of $\SS$ generated by bad arcs, is a subalgebra of $\extX_+(\surface;\quasi)$. As $\SS$ is spanned over $\cD(\fS)$ by $\Bo^\rd$ (by Corollary \ref{r.badbasis}), and bad arcs $q$-commute with monomials $x^\bk$ (by Lemma \ref{r.badarc}), we can assume that $\al\in \Bo^\rd$.

Assume that $\al$ has a negative state on some boundary edge $e$. Consider the product $X_e \al$, whose diagram is depicted in Figure~\ref{fig:bad2}, where we also depict $Y_+$ and $Y_-$ which are respectively the result of positive and negative resolution of the diagram of $X_e\al$ at the lowest crossing.

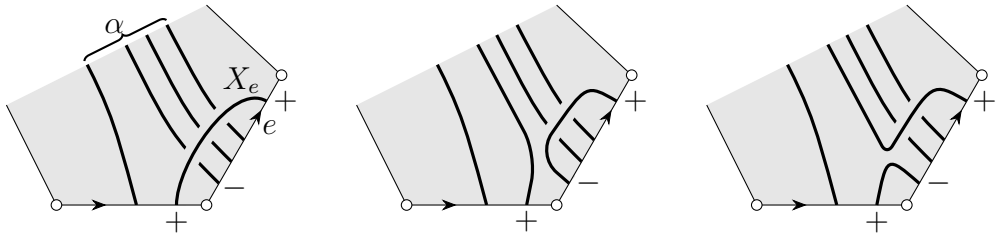
\begin{figure}[h]
\centering
\begin{tikzpicture}[scale=1.33]
\fill[gray!20!white] (-2,1)--(-1.5,0)--(0,0)--(60:1.5)--(0,2);
\begin{scope}
\clip (-2,1)--(-1.5,0)--(0,0)--(60:1.5)--(0,2);
\draw[very thick] (-0.7,0)..controls ++(105:1)..(-1.5,2)
	(60:0.25)..controls ++(135:0.7)..(-1,2)
	(60:0.5)..controls ++(135:0.6)..(-0.75,2)
	(60:0.75)..controls ++(135:0.5)..(-0.5,2);
\end{scope}
\draw[knot,very thick,background color=gray!20!white,knot gap=6,inner sep=2pt]
	(60:1.2)node[right]{$+$}..controls ++(150:0.4) and (-0.3,0.4)
	..(-0.3,0)node[below]{$+$};
\draw (-2,1)--(-1.5,0)--(0,0)--(60:1.5)--(0,2);
\path[tips,-{Stealth[length=0.2cm]}] (-1.5,0)--(-1,0);
\path[tips,-{Stealth[length=0.2cm]}] (0,0)--(60:1.1);
\draw[fill=white] (0,0)circle(1.5pt) (-1.5,0)circle(1.5pt) (60:1.5)circle(1.5pt);
\draw (60:1.5)++(-0.4,-0.05)node{$X_e$} (60:0.9)node[right]{$e$}
	(60:0.2)node[inner sep=2pt,right]{$-$};
\draw[thick,decorate,decoration=brace] (-1.22,1.44)--(-0.4,1.85);
\draw (-0.85,1.65)++(-0.07,0.14)node{$\alpha$};

\begin{scope}[xshift=3.5cm]
\fill[gray!20!white] (-2,1)--(-1.5,0)--(0,0)--(60:1.5)--(0,2);
\begin{scope}
\clip (-2,1)--(-1.5,0)--(0,0)--(60:1.5)--(0,2);
\draw[very thick] (-0.7,0)..controls ++(105:1)..(-1.5,2)
	(60:0.5)..controls ++(135:0.6)..(-0.75,2)
	(60:0.75)..controls ++(135:0.5)..(-0.5,2);
\draw[very thick,rounded corners=10pt] (-0.3,0)--(-0.2,0.5)--(-1,2);
\end{scope}
\draw[knot,very thick,background color=gray!20!white,knot gap=6,inner sep=2pt,rounded corners=10pt]
	(60:1.2)node[right]{$+$}..controls ++(150:0.35)..(-0.2,0.5)--
	(60:0.25)node[right]{$-$};
\draw[inner sep=2pt] (-0.3,0)node[below]{$+$};
\draw (-2,1)--(-1.5,0)--(0,0)--(60:1.5)--(0,2);
\path[tips,-{Stealth[length=0.2cm]}] (-1.5,0)--(-1,0);
\path[tips,-{Stealth[length=0.2cm]}] (0,0)--(60:1.1);
\draw[fill=white] (0,0)circle(1.5pt) (-1.5,0)circle(1.5pt) (60:1.5)circle(1.5pt);
\end{scope}

\begin{scope}[xshift=7cm]
\fill[gray!20!white] (-2,1)--(-1.5,0)--(0,0)--(60:1.5)--(0,2);
\begin{scope}
\clip (-2,1)--(-1.5,0)--(0,0)--(60:1.5)--(0,2);
\draw[very thick] (-0.7,0)..controls ++(105:1)..(-1.5,2)
	(60:0.5)..controls ++(135:0.6)..(-0.75,2)
	(60:0.75)..controls ++(135:0.5)..(-0.5,2);
\draw[knot,very thick,background color=gray!20!white,knot gap=6,rounded corners] (60:1.2)..controls ++(150:0.4)..(-0.2,0.5)--(-1,2);
\end{scope}
\draw[very thick,inner sep=2pt]
	(-0.3,0)node[below]{$+$}..controls (-0.2,0.5)..(60:0.25)node[right]{$-$};
\draw[inner sep=2pt] (60:1.2)node[right]{$+$};
\draw (-2,1)--(-1.5,0)--(0,0)--(60:1.5)--(0,2);
\path[tips,-{Stealth[length=0.2cm]}] (-1.5,0)--(-1,0);
\path[tips,-{Stealth[length=0.2cm]}] (0,0)--(60:1.1);
\draw[fill=white] (0,0)circle(1.5pt) (-1.5,0)circle(1.5pt) (60:1.5)circle(1.5pt);
\end{scope}
\end{tikzpicture}
\caption{Diagrams of $X_e \al$ and $Y_+, Y_-$.}\label{fig:bad2}
\end{figure}

By the skein relation $X_e \al= q Y_++ q^{-1} Y_-$. Since $Y_-$ has a bad arc, Lemma~\ref{r.badarc} shows that $Y_-$ is an $\cD(\fS)$-linear combination of elements in $\Bo^\rd$ with less than $d(\al)$ negative states. The diagram $Y_+$ contains an arc whose two endpoints are in $e$, and by resolving all the crossing on this arc and using the defining relation \eqref{eq.arcs}, we get that $Y_+$ is also an $\cD(\fS)$-linear combination of elements in $\Bo^\rd$ with less than $d(\al)$ negative states. Thus by induction, the lemma holds for $\al$. This completes the proof of Lemma~\ref{lemma-Laurent} and the theorem.
\end{proof}

\begin{remark} \label{r.bad5}
The induction proof of Lemma~\ref{lemma-Laurent} shows that for every $\al\in \SS$ there is a monomial $m$ in the variables $X_e, e\in \Ed$ such that $m \al\in\stateS_\bad^+(\fS)$, the $\RV$-subalgebra of $\SS$ generated by bad arcs and $\cS^+(\fS)$.
\end{remark}

\section{Reduced skein algebra and quantum cluster algebra}\label{sec.red}

\subsection{Reduced skein algebra} Suppose $\fS=\bfS \setminus \cP$ is a punctured bordered surface with non-empty boundary. Let $I^\bad$ be the left ideal generated by all bad arcs defined in Subsection \ref{sec.bad}. By Lemma~\ref{r.badarc}, the ideal $I^\bad$ is a two-sided ideal. The quotient $\cSd(\surface)=\stateS(\surface)/I^\bad$ is the \emph{reduced skein algebra} introduced in \cite{CL}, where it is shown that $I^\bad$ is the kernel of Bonahon-Wong'quantum trace map, see Section \ref{sec-shear}.


Recall the basis of $\stateS(\surface)$ is given by $B(\surface)$, the set of isotopy classes of increasingly stated, positively-ordered simple $\partial\surface$-tangle diagrams. Let $B^\mathrm{rd}(\surface)\subset B(\surface)$ be the subset of diagrams which contain no bad arcs. 

\begin{theorem}[Theorems 7.1 and 7.12 in \cite{CL}] Suppose $\fS$ is a punctured bordered surface with non-empty boundary. The algebra $\reduceS(\surface)$ is a domain, and
$B^\mathrm{rd}(\surface)$ is a free $\mathcal{R}$-basis of $\reduceS(\surface)$.
\end{theorem}

\def\cA{\mathcal A}
\subsection{Quantum cluster algebra}
Recall that the Muller skein algebra $\cS^+(\fS)\subset \SS$ is the $\cR$-subalgebra generated by $\pfS$-tangle diagrams with positive states. For every ideal arc $e$ in $\fS$ we defined an element $X_e$, which is the result of moving endpoints of $e$ to the left and equipping them with positive states at both endpoints, see Subsection \ref{sec.bad}.

Let $\fM$ denote the multiplicative subset generated by $X_e$ where $e$ runs over the set $\quasi_\partial$ of boundary edges. Muller showed in \cite{Muller} that $\fM$ is a two sided Ore set in $\mullerS(\surface)$. Let $\cA(\fS)$ be the $\cR$-subalgebra of the Ore localization $\mullerS(\surface)\fM^{-1}$ generated by all $X_a$, where $a$ can be any ideal arc, and the inverses of $X_e, e\in \Ed$. When $\fS$ has no interior points, Muller \cite{Muller} showed that $\cA(\fS)$ is a quantum cluster algebra (in the sense of \cite{BZ}) which quantizes the classical cluster algebra associated with $(\overline{\surface},\marked)$ defined in \cite{GSV,FG}. If in addition each connected component of $\fS$ is triangulable and has at least two puncture, then Muller showed that $\cA(\fS)= \mullerS(\surface)\fM^{-1}$.

\subsection{Relation between reduced skein algebra and quantum cluster algebra}

The set $B^+(\fS)\subset B(\fS)$ consisting of diagrams with positive states is an $\cR$-basis of $\cS^+(\fS)$. As $B^+(\fS) \subset B^\rd(\fS)$, the inclusion $\mullerS(\surface)\hookrightarrow\stateS(\surface)$ descends to an injective map
\begin{equation}\label{eq.incl1}
\varkappa: \mullerS(\surface)\hookrightarrow\reduceS(\surface).
\end{equation}

While $\mullerS(\surface)$ and its localization $\mullerS(\surface)\mathfrak{M}^{-1}$ are defined using only positive states, the reduced stated skein algebra $\reduceS(\surface)$ uses both positive and negative states. Hence it is a surprise that we can have the following, which demonstrates the ubiquity of the stated skein algebra and gives another perspective for the quantum cluster algebra $A(\fS)$.

\def\tk{\tilde \varkappa}

\begin{theorem}\label{thm.reduced-cluster}
Suppose $\fS$ is a connected punctured bordered surface with non-empty boundary. The embedding $\varkappa$ extends uniquely to an $\cR$-algebra isomorphism
$$\tk: \mullerS(\surface)\fM^{-1} \to \reduceS(\surface).$$

Consequently, when $\fS$ has at least two punctures and no interior punctures, the quantum cluster algebra $\cA(\fS)$ is naturally isomorphic to the reduced skein algebra $\reduceS(\surface)$.
\end{theorem}

\begin{proof} 
A crucial fact \cite[Proposition 7.4]{CL} is that the elements $X_e, e\in \Ed$ are invertible in $\reduceS(\surface)$. The inverse of $X_e$ is given by the same arc, only with negative states on both endpoints. The universality of the Ore extension shows that there is a unique $\cR$-algebra extension
\[\tk: \mullerS(\surface)\fM^{-1} \to \reduceS(\surface).\]

To show $\tk$ is surjective, for any $\alpha\in\reduceS(\surface)$, choose a lift $\alpha'\in\stateS(\surface)$. There exists $m\in\mathfrak{M}$ such that $m\alpha'\in\stateS_\bad^+(\surface)$ by Remark~\ref{r.bad5}, where $\stateS_\bad^+(\fS)$ is the $\RV$-subalgebra of $\SS$ generated by bad arcs and $\cS^+(\fS)$. Passing to the reduced algebra, $m\alpha\in\mullerS(\surface)$. Thus $\alpha=m^{-1}(m\alpha)$ is in the image of $\tk$.

Identify $\mullerS(\surface)$ with a subset of $\reduceS(\surface)$. Since every $m\in \fM$ is invertible in $\reduceS(\surface)$ and the map $\tk$ is surjective, $\reduceS(\surface)$ is a left ring of fractions of $\mullerS(\surface)$ with respect to $\fM$. By the uniqueness of left ring of fractions, see \cite[Corollary 6.4]{GW}, the map $\tk$ is an isomorphism.
\end{proof}

\section{Embedding into quantum torus, shear coordinate version}
\label{sec-shear}

In this section, $\surface$ is a punctured bordered surface with at least one puncture. Unlike the case considered in Section \ref{sec.length}, we do not require $\fS$ to have non-empty boundary.

We will show that there is an algebra embedding of $\SS$ into a quantum torus
\[\phi:\stateS(\surface)\hookrightarrow\mathbb{T}(\bar{Q}).\]
The matrix $\bar Q$ depends on an ideal triangulation of $\fS$.

When $\fS$ does not have boundary, then $\SS$ is the ordinary skein algebra and the map $\phi$ was constructed by Bonahon and Wong \cite{BW1}. In this case, when $q=1$ the image $\phi(\al)$ of a simple closed curve $\al$ expresses the trace of $\al$ as a Laurent polynomial in the shear coordinates of the enhanced Teichm\"uller space.

Besides, we will show that $\CS$ is Noetherian, orderly finitely generated, and calculate the Gelfand-Kirillov dimension of $\CS$.

\subsection{Ideal triangulation and Chekhov-Fock algebra}

A punctured bordered surface is \emph{triangulable} if it has at least one puncture, and it is not one of the following: a monogon, a bigon, and the sphere with 1 or 2 punctures. In this section we assume that $\fS$ is triangulable.

An \emph{ideal triangulation} of $\fS$ is a maximal collection $\Delta$ of non-trivial ideal arcs which are pairwise disjoint and pairwise non-isotopic.
Every boundary edge is isotopic to an element of $\Delta$, and distinct boundary edges are not isotopic. After an isotopy, we can assume that $\Delta$ contains the set $\Dd$ of all boundary edges. Of course $\Dd$ is the same as $\Ed$ of Section~\ref{sec.length}, but since we use the notation $\Delta$ for the set of edges of the triangulation, we use $\Dd$ for the set of boundary edges. Let $\Do=\Delta\setminus \Dd$ be the set of all \emph{interior edges} of the triangulation.

By splitting $\fS$ along all interior edges in $\Do$, we get a collection $\cF$ of ideal triangles with the projection
$p: \bigsqcup_{\fT\in\cF}\fT \onto \fS$.
If $a,b,c$ are the edges of a triangle $\fT$ in $\cF$, then we call $(p(a), p(b), p(c))$ a \emph{triangular triple}. Two of the edges may coincide, in which case we call $\fT$ a \emph{self-folded triangle}, and the repeated edge is called a \emph{self-folded edge}. $p(\fT)$ is a punctured monogon in this case.

\newcommand{\facedef}[2]{
\begin{tikzpicture}[scale=0.8,baseline=0.28cm]
\draw[fill=gray!20!white] (0,1)--(0.6,0)--(1,1);
\draw[inner sep=0pt] (0.1,0.5)node{\vphantom{$b$}#1} (1,0.5)node{\vphantom{$b$}#2};
\draw[fill=white] (0.6,0)circle(2pt);
\end{tikzpicture}
}

\def\hDd{{\hat \Delta_\partial}}

The Thurston form of $\Delta$ is the antisymmetric function $Q_{\Delta}: \Delta\times \Delta\to \BZ$ defined by
\begin{equation}\label{eq.Q}
Q_\Delta(a,b) = \#\left( \facedef{$b$}{$a$} \right)- \#\left( \facedef{$a$}{$b$} \right).
\end{equation}
Here each shaded part is a corner of an ideal triangle. Thus, the right hand side of \eqref{eq.Q} is the number of corners where $b$ is counterclockwise to $a$ minus the number of corners where $a$ is clockwise to $b$, as viewed from the puncture. The Thurston form is related to the Weil-Petersson Poisson structure of the enhanced (or holed) Teichm\"uller space in shear coordinates, see e.g. \cite{BW1,FG2}. The matrix $Q_\Delta$ is called the face matrix in \cite{Le:qtrace}.

The Chekhov-Fock algebra $\cY(\fS;\Delta)$ is the quantum torus associated to
$Q_\Delta$
\begin{equation}
\cY(\fS;\Delta):= \cR\la y_a^{\pm1}, a \in \Delta \ra / ( y_ a y_b = q^{Q_\Delta(a,b) } y_b y_a).
\notag
\end{equation}

The set $\{ y^{\bk} \mid \bk \in \BZ^\Delta\}$ is a free $\cR$-basis of $\cY(\fS;\Delta)$. A map $\bk:\Delta\to \BZ$ is \emph{balanced} if $\bk(a) + \bk(b) + \bk(c)$ is even whenever $a,b,c$ are a triangular triple. Let $\Ybl(\fS;\Delta)$ be the monomial subalgebra of $\cY(\fS;\Delta)$ spanned by all $y^\bk$ with balanced $\bk$. Then $\Ybl(\fS;\Delta)$ is itself a quantum torus. Let $\Fr(\Ybl(\fS;\Delta))$ be the division algebra of fractions of $\cY(\fS;\Delta)$.

Hiatt \cite{Hiatt} and Bonahon and Wong \cite{BW1}, extending the work of Chekhov-Fock \cite{CF} and Liu \cite{Liu}, showed that given another ideal triangulation $\Delta'$, there is a natural coordinate change isomorphism
\[\Theta_{\Delta\Delta'}: \Fr(\Ybl(\fS;\Delta')) \to \Fr(\Ybl(\fS;\Delta))\] such that $\Theta_{\Delta\Delta}=\id$ and $\Theta_{\Delta\Delta'}\circ \Theta_{\Delta'\Delta''} = \Theta_{\Delta\Delta''}$.

Let us discuss a criterion for a monomial to be balanced. Suppose $\al$ is a $\pfS$-tangle diagram. For an ideal arc $a$ the \emph{geometric intersection number} $I(a,\al)$ is the minimum among all $\abs{\al' \cap a}$, where $\al'$ is isotopic to $\al$. We say $\al$ is \emph{taut} with respect to a finite collection $X$ of disjoint ideal arcs in $\fS$ if $\abs{\al \cap a}= I(\al,a)$ for all $a\in X$.

Define $\bn_\al: \Delta\to \BN$ by $\bn_\al(a) = I(\al,a)$ for $ a\in \Delta$. It is easy to show that $\bn_\al$ is balanced. Besides, any element in $2\BZ^\Delta$ is clearly balanced. The following characterization of balanced elements, which is easy to prove, is essentially contained in \cite[Section 2]{BW1}.

\begin{lemma}\label{r.bl1}
An element $\bk\in \BZ^\Delta$ is balanced if and only there is a $\pfS$-tangle diagram $\al$ such that $\bk - \bn_\al \in 2 \BZ^\Delta$.
\end{lemma}

\subsection{Bonahon-Wong quantum trace and its descendants}

Recall that the Bonahon-Wong version $\tSS$ is the skein algebra defined just like $\SS$, where one uses only the skein relation \eqref{eq.skein} and the trivial knot relation \eqref{eq.loop}. Thus our $\SS$ is the quotient of $\tSS$, factored out by the other two relations \eqref{eq.arcs} and \eqref{eq.order}.

The quantum trace map Bonahon and Wong is an algebra homomorphism
\[\tr^\BW_\Delta: \tSS \to \Ybl(\fS;\Delta),\]
compatible with the coordinate change isomorphisms: for another ideal triangulation $\Delta'$,
\begin{equation}\label{eq.comp}
\tr^\BW_\Delta = \Theta_{\Delta \Delta'}\circ \tr^\BW_{\Delta'}.
\end{equation}

When $\fS$ has no boundary, the map $\tr^\BW_\Delta$ is injective \cite{BW1}. However, if $\pfS\neq \emptyset$, then $\tr^\BW_\Delta$ is not injective. It is proved in \cite{Le:TDEC} that $\tr^\BW_\Delta$ descends to an algebra homomorphism
\[\tr_\Delta: \cS(\fS ) \to \Ybl(\fS;\Delta).\]
The map $\tr_\Delta$ is still not injective. In \cite{CL} it is proved that the kernel of $\tr_\Delta$ is the ideal $I^\bad$ generated by bad arcs. Hence $\tr^\BW_\Delta$ further descends to an algebra \emph{embedding}
\begin{equation}
\tr_\Delta: \cS^\rd(\fS ) \embed \Ybl(\fS;\Delta).
\end{equation}

\no{
\begin{remark}
When $q=1$ and $\fS$ has no boundary, the ideal triangulation $\Delta$ identifies the enhanced Teichm\"uller space of $\fS$ with a Euclidean space via shear coordinates, and the image $\tr_\Delta(\al)$ of a simple closed curve $\al$ expresses the trace of $\rho(\al)$, where $\rho$ is the $PSL_2(\BR)$-representation of the fundamental group $\pi_1(\fS)$ corresponding to the hyperbolic structure, as a Laurent polynomial in the shear coordinates. See \cite{BW1} for the precise statement and details. For this reason $\tr^\BW_\Delta$ is called the quantum trace.
\end{remark}
}

\subsection{Extended Chekhov-Fock algebra}\label{sec.ext_CF}

To embed $\SS$ into a quantum torus, we extend the balanced Chekhov-Fock algebra $\Ybl(\fS;\Delta)$ to a bigger quantum torus.

Let $\hDd=\{ \hat e \mid e\in \Dd\}$ be another copy of $\Dd$ and $\bD= \Delta\sqcup \hDd$. Extend the Thurston form $Q_\Delta$ to an antisymmetric function $\bQ_\Delta:\bD \times \bD \to \BZ$ so that the values of $\bQ_\Delta$ on the extension set $(\bD \times \bD) \setminus (\Delta\times \Delta)$ are 0 except
\[\bQ_\Delta(\hat e , e)=1=-\bQ_\Delta(e,\hat e) \text{ for all } e\in \Dd.\]
The geometric origin of $\bQ_\Delta$ will be given in Subsection~\ref{sec.coord}. Consider the quantum torus
\begin{equation}\label{eq.TQbar}
\bT(\bQ_\Delta)= \cR\la z_a^{\pm1}, a \in \bD \ra / ( z_ a z_b = q^{\bQ_\Delta(a,b) } z_b z_a).
\end{equation}
Identify $\cY(\fS;\Delta)$ with a subalgebra of $\bT(\bQ_\Delta)$ via the embedding $y_a \mapsto z_a, a\in \Delta$. Let $\bY^\bal(\fS;\Delta)$ be the subalgebra of $\bT(\bQ)$ generated by $\Ybl(\fS;\Delta)$ and all the $z_\he^{\pm 2}, e\in \Dd$.
We have an $\cR$-linear map $\pr: \bY^\bal(\fS;\Delta) \onto \cY^\bal(\fS;\Delta)$ defined by
\begin{equation}\label{eq.pr}
\pr(z^\bk) = \begin{cases}
z^\bk, &\text{if }\bk(\he)= 0\text{ for all } e\in \Dd,\\
0, &\text{otherwise},
\end{cases}
\end{equation}
which is identity on $\Ybl(\fS;\Delta)$.

Let $\cZ(\fS;\Delta)$ be the subalgebra of $\bT(\bQ)$ generated by $\balancedY(\fS;\Delta)$ and $z_\he^2, e\in \Dd$. Then the restriction $\pr:\cZ(\fS;\Delta)\to \balancedY(\fS;\Delta)$ is an $\cR$-algebra homomorphism.


For a stated simple $\pfS$-tangle diagram $\al$ with state $s$, define $\bbn_\al: \bD \to \BN$ by:
\begin{equation}\label{eq.bbn}
\begin{cases} \bbn_\al(a):= \bn_\al(a)=I(\al,a), & a\in \Delta,\\
\bbn_\al(\he):=\abs{\al\cap e} - \sum_{x\in (\al\cap e)} s(x),& e\in\Delta_\partial.
\end{cases}
\end{equation}
From Lemma~\ref{r.bl1} one has the following characterization of $\bYbl(\fS;\Delta)$ and $\cZ(\fS;\Delta)$.

\begin{lemma}\label{r.bl}
Suppose $\bk\in \BZ^\bD$. Then $z^\bk\in \bYbl(\fS;\Delta)$ if and only if there is a stated simple $\pfS$-tangle diagram $\al$ such that $\bk -\bbn_\al \in 2 \BZ^\bD$. Furthermore $z^\bk\in \cZ(\fS;\Delta)$ if and only if in addition $\bk(\he) \ge 0$.
\end{lemma}

\subsection{Parameterization of basis $B(\fS)$}

By Theorem \ref{thm.basis}, the set $B(\fS)$ of isotopy classes of increasingly stated, positively ordered, simple $\pfS$-tangle diagrams is an $\cR$-basis of $\CS$.

\begin{proposition}\label{r.para}
Suppose $\Delta$ is a triangulation of a triangulable punctured bordered surface $\fS$. The map $B(\fS) \to \BZ^\bD$ given by $\al \to \bbn_\al$ is injective, and its image is the submonoid $\Lambda_\Delta\subset \BZ^\bD$ consisting of $\bn\in \BN^\bD$ such that
\begin{enumerate}
\item for any triangular triple $a,b,c\in \Delta$, $\bn(a) + \bn(b) + \bn(c) \in 2\BN$ and $ \bn(a) \le \bn(b)+ \bn(c)$,
\item for $e\in \Dd$, $\bn(\he) \in 2\BN$ and $\bn(\he) \le 2 \bn(e)$.
\end{enumerate}
Moreover, the rank of $\Lambda_\Delta$ 
is $r(\surface):=\abs{\bar{\Delta}}$.
\end{proposition}

\begin{remark}
Let $\chi(\fS)$ be the Euler characteristic and $\abs{\Pd}$ the number of boundary punctures. An easy counting argument shows that
\begin{equation}\label{eqn-r-def}
r(\fS) = 3 \abs{\Pd}-3 \chi(\fS).
\end{equation}

\end{remark}

\begin{proof}
If $\al\in B(\fS)$ it is easy to show that $\bbn_\al\in \Lambda_\Delta$.
In the opposite direction, it is well known \cite{Matveev} that (a) implies that there is a unique (non-stated) simple $\pfS$-tangle diagram $\al$ such that $\bn_\al(a) = \bn(a)$ for $a\in \Delta$. The number of $+$ and $-$ states of $\al\cap e$ are also uniquely determined by $\bn(\he)$ given (b), and the states are completely fixed for increasingly stated diagrams. Thus for every $\bn\in \Lambda_\Delta$, there is a unique $\al\in B(\fS)$ such that $\bbn_\al= \bn$.

To compute the rank, we show that the group generated by $\Lambda_\Delta$ contains $(2\mathbb{Z})^{\bar{\Delta}}$. Let $\mathbf{2}\in\mathbb{Z}^{\bar{\Delta}}$ be the constant map $2$, and $\mathbf{d}_a\in\mathbb{Z}^{\bar{\Delta}}$ be the indicator function on $a\in\bar{\Delta}$. It is clear that $\mathbf{2}$ and $\mathbf{2}+2\mathbf{d}_a$ are in $\Lambda$ for all $a\in\bar{\Delta}$. Thus the difference $2\mathbf{d}_a$ is in the group generated by $\Lambda_\Delta$, and they span $(2\mathbb{Z})^{\bar{\Delta}}$.
\end{proof}

\subsection{Filtrations} \label{sec.filtr}

Define a group homomorphism $\dego:\mathbb{Z}^{\bar{\Delta}}\to\mathbb{Z}$ by
\[\dego(\mathbf{k})=\sum_{e\in\Delta}\mathbf{k}(e).\]
Note that $\mathbf{k}(\hat e)$, $e\in \Dd$, are not taken into account. Then $\dego$ induces a $\mathbb{Z}$-grading of $\mathbb{T}(\bar{Q}_\Delta)$ from \eqref{eq.grad}, and it further induces a filtration as follows. For $k\in\mathbb{N}$, let $F_k(\bT(\bQ_\Delta))\subset \bT(\bQ_\Delta)$ be the $\cR$-submodule spanned by $z^\bk$ such that $\dego(\mathbf{k})\le k$. Then $\{F_k(\bT(\bQ_\Delta))\}$ is a filtration of $\bT(\bQ_\Delta)$ compatible with the algebra structure. On a subalgebra $A$, there is an induced filtration by letting $F_k(A):= F_k(\bT(\bQ_\Delta)) \cap A$.

For a stated simple $\pfS$-tangle diagram $\al$ define $\dego(\al)\eqdef \dego(\bbn_\al)=
\sum_{a\in\Delta} \bn_\al(a)$. Note that $\bbn_\al(\he)$, $e\in \Dd$, are not taken into account. For $k\in \BN$, let $F_k(\CS)\subset \CS$ be the $\cR$-submodule spanned by stated simple $\pfS$-tangle diagrams $\al$ such that $\dego(\al)\le k$. Then $\{F_k(\CS)\}$ is an $\BN$-filtration of $\cS(\fS)$, which is compatible with the algebra structure. The set $B_k(\fS):=\{\al\in B(\fS) \mid \dego(\al) \le k\}$ is a free $\cR$-basis of $F_k(\CS)$.

\subsection{Quantum trace map, GK dimension, and orderly finite generation}

\begin{theorem}\label{thm.embed3}
Suppose $\fS$ is a triangulable punctured bordered surface
and $\Delta$ is an ideal triangulation of $\fS$.

\begin{enumerate}
\item There is an algebra embedding
\begin{equation}\label{eq.qtr2}
\phi_\Delta: \cS(\fS) \embed \cZ(\fS;\Delta )
\end{equation}
which is a lift of the quantum trace $\tr_\Delta$, i.e., the map $\tr_\Delta$ is the composition
\[\cS(\fS) \xrightarrow {\phi_\Delta} \cZ(\fS;\Delta ) \xrightarrow {\pr} \Ybl(\fS;\Delta ).\]
In addition, if $\cR=\BZ[q^{\pm 1/2}]$, then $\phi_\Delta$ is reflection invariant.

\item The homomorphism $\phi_\Delta$ is compatible with the filtrations defined in Subsection \ref{sec.filtr}. Moreover, if $\alpha$ is a stated simple $\pfS$-tangle diagram, then there is $t(\al)\in \frac12 \BZ$ such that
\begin{equation}\label{eqn-tr-highest}
\phi_\Delta(\alpha)=q^{t(\al)} z^{\bbn_\al} \mod F_{\dego(\al) -1}(\cZ(\fS;\Delta)).
\end{equation}

\item For another ideal triangulation $\Delta'$, the coordinate change $\Theta_{\Delta\Delta'}:\Fr(\Ybl(\fS;\Delta'))\to \Fr(\Ybl(\fS;\Delta))$ extends to a unique algebra isomorphism
$$ \bTheta_{\Delta\Delta'}:\Fr(\bYbl(\fS;\Delta'))\to \Fr(\bYbl(\fS;\Delta))$$
given on extra variables $z_{\hat e} $ by $\bTheta(z_{\hat e}^2)= z_{\hat e}^2$ for all $e\in \Dd$. For yet another ideal triangulation $\Delta''$, we have $\bTheta_{\Delta''\Delta'} \circ \bTheta_{\Delta'\Delta} = \bTheta_{\Delta''\Delta}$ and $\bTheta_{\Delta \Delta}=\id$.

The embedding $\phi_\Delta$ is compatible with coordinate change, i.e.,
\begin{equation}\label{eq.comp1}
\bTheta_{\Delta' \Delta} \circ \phi_\Delta= \phi_{\Delta'}.
\end{equation}
\end{enumerate}
\end{theorem}

As consequences, we have the following two statements.

\begin{theorem}\label{r.degen}
With the assumption of Theorem \ref{thm.embed3}, the associated graded algebra of $\CS$ with respect to the filtration $F_k(\CS)$ is isomorphic to the monomial subalgebra $A(\bQ;\Lambda_\Delta)$ of the quantum torus $\bT(\bQ)$, where the submonoid $\Lambda_\Delta$ is defined in Proposition~\ref{r.para}.
\end{theorem}


\begin{theorem}\label{thm-noether1}
Let $\fS$ be a punctured bordered surface. 
\begin{enumerate}
\item As an $\cR$-algebra, $\cS(\fS)$ is orderly finitely generated by arcs and loops. This means, there are one-component simple stated $\pfS$-tangle diagrams $\alpha_1, \dots, \alpha_n$ such that the set $\{ \alpha_1^ {k_1} \dots \alpha_n^{k_n} \mid k_i \in \BN\}$ spans $\cS(\fS)$ over $\cR$.
\item $\cS(\fS)$ is a Noetherian domain.
\item If $\fS$ is triangulable, then
the GK dimension of $\cS(\fS)$ is $r(\fS)$ defined by \eqref{eqn-r-def}.
\end{enumerate}
\end{theorem}


Part (c) of Theorem~\ref{thm-noether1} means the domain and the target space of the embedding \eqref{eq.qtr2} have the same GK dimension. Hence $\SS$ cannot embed into a quantum torus of less dimension.


\no{
\begin{remark}
When $\pfS=\emptyset$, then $\SS$ is the ordinary skein algebra, and most of Theorem \ref{thm-noether1} was known: the finite generation (without order) was proved in \cite{Bullock}, the orderly finite generation was proved in \cite{AF}, and the Noetherian domain property was established in \cite{PS2}.
\end{remark}
}

\begin{remark}
When $\fS$ is the closed surface of genus $g$, the method of this paper does not apply as there is no ideal triangulation. In this case one can still show that GK dimension of $\SS$ is still $-3\chi(\Sigma_g)=6g-6$, as defined in \eqref{eqn-r-def}. For a proof see \cite{KL}.
\end{remark}

\begin{remark}\label{rmk-qpower}
In part (b) of Theorem \ref{thm.embed3}, the exponent $t(\al)$ can be easily determined using the reflection invariance property. In particular, if the diagram has at most one endpoint on each boundary edge, then $t(\al)=0$ since the diagram is already reflection invariant.
\end{remark}

In the remaining of this section we give proofs of Theorems~\ref{thm.embed3}--\ref{thm-noether1}. The strategy is to embed $\fS$ into a bigger surface $\fS^\ast$ and utilize the quantum trace map of $\fS^\ast$. First we recall the quantum trace map (via splitting) and some of its properties in Subsections~\ref{sec.triangle}--\ref{sec.BW}.

\def\bm{{\mathbf m}}

\subsection{Ideal triangle} \label{sec.triangle}

\begin{figure}[h]
\centering
\begin{tikzpicture}
\draw[fill=gray!20!white] (-60:2.5) to node[right,pos=0.35]{$b$} (0,0)
	to node[left,pos=0.65]{$c$} (-120:2.5)
	to node[below]{$a$} cycle;
\draw[fill=white] (-60:2.5)circle(2pt) (0,0)circle(2pt)
	(-120:2.5)circle(2pt);
\draw[very thick] (-120:0.75)node[left]{$\mu$} to[out=-30,in=-150] node[below]{$\alpha$} (-60:0.75)node[right]{$\nu$};
\end{tikzpicture}
\hskip2cm
\begin{tikzpicture}
\draw[fill=gray!20!white] (-60:2.5) to node[right]{$b$} (0,0)
	to node[left]{$c$} (-120:2.5) to node[below]{$a$} cycle;
\draw[fill=gray!20!white]
	(-60:2.5) to[out=140,in=-80] node[left]{$b'$} (0,0)
	to[out=-100,in=40] node[right,pos=0.65,inner sep=2pt]{$c'$} (-120:2.5)
	to[out=20,in=160] node[above,inner sep=1pt]{$a'$} cycle;
\draw[fill=white] (-60:2.5)circle(2pt) (0,0)circle(2pt)
	(-120:2.5)circle(2pt);
\end{tikzpicture}
\caption{Left: Ideal triangles with arc $\al(\mu,\nu)$. Right: ideal arcs $a',b',c'$.}\label{fig:tria1}
\end{figure}
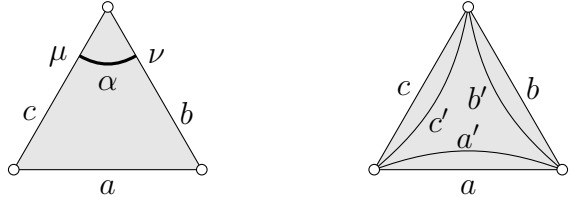

Let $\fT$ be an ideal triangle, with boundary edges $a,b,c$ in counterclockwise order as in Figure~\ref{fig:tria1}.
Then $\Delta=\{a,b,c\}$ is the unique ideal triangulation of $\fT$, and $\cY(\fT)=\cY(\fT;\Delta)$ is the quantum torus
\begin{equation}
\cY(\fT) : = \cR \la y_a^{\pm 1},y_b ^{\pm 1}, y_c^{\pm 1} \ra /( qy_a y_b = y_by_a, qy_by_c = y_cy_b, qy_cy_a = y_ay_c).
\end{equation}

\def\ST{\cS(\fT)}
\def\YT{\cY(\fT)}
\def\pfT{\partial \fT}
\def\Coeff{\operatorname{Coeff}}

The ideal triangle $\fT$ has an order $3$ rotational symmetry. The algebra $\ST$ has $12$ generators given by $\alpha(\mu,\nu)$, $\mu,\nu\in\{\pm\}$, as in Figure~\ref{fig:tria1}, and their images under the rotational symmetries. The quantum trace map $\tr_\fT:\stateS(\mathfrak{T})\to\fockY(\fT)$ is the reflection invariant algebra homomorphism, equivariant under rotation, defined by
\begin{equation}\label{eq.triangle}
\tr(\alpha(\mu,\nu))=\begin{cases}
0,&(\mu,\nu)=(-,+),\\
[y_c^\mu y_b^\nu],&\text{otherwise}.
\end{cases}
\end{equation}
Recall $\pm$ are identified with $\pm1$. Note that $\tr(\alpha(\mu,\nu))=0$ exactly when $\alpha(\mu,\nu)$ is a bad arc.

\begin{lemma}\label{r.tri5}
Suppose $\al$ is a stated simple $\pfT$-tangle diagram with state $s:\partial\alpha\to\{\pm\}$. For $e\in\{a,b,c\}$, let $\bk(e)= \sum_{x\in e\cap \al} s(x)$. Then
\begin{equation}\label{eq.tri1}
\tr_\fT(\al) = \Coeff(\al) y^{\bk}, \quad \Coeff(\al) \in \cR.
\end{equation}
In addition, if all the states of $\al$ are positive, then $\Coeff(\al)= q^t$ for some $t\in \frac12 \BZ$.
\end{lemma}

\begin{proof}
For $e\in \{a,b,c\}$, choose an ideal arc $e'$ isotopic to $e$. Splitting along $e'$ creates a bigon $\fB_e$ as in Figure~\ref{fig:tria1}. We assume $\al$ is taut with respect to $a', b', c'$. Let $\fT'$ be the triangle bounded by $a',b',c'$. The obvious diffeomorphism from $\fT'$ to $\fT$ gives an identification $\ST\equiv \cS(\fT')$. Choose a linear order of connected components of $\al\cap \fT'$. Isotope the height of $\alpha$ so that each connected component of $\alpha\cap\mathfrak{T}'$ has constant height, and the heights satisfy the linear order. This also induces a linear order on $\al\cap e'$ for each $e\in \{a,b,c\}$. Applying the splitting homomorphism on $\fT$ along $a',b',c'$ and using the counit on the bigons, we get that the composition
\[\ST \to \cS(\fT') \ot \left(\cS(\fB_a) \ot \cS(\fB_b) \ot \cS(\fB_c) \right) \xrightarrow{\id \ot \ve} \ST\]
is the identity. It follows that
\begin{equation}
\tr_\fT(\al) = \sum_{s'} \{\tr_{\fT'}(\al\cap \fT',s')\} \prod_{e\in \{a,b,c\}} \left[ \ve ( \al\cap \fB_e,s' ) \right],
\end{equation}
where $s'$ runs over the set $s': \{a',b',c'\}\to \{\pm \}$, and $(\al\cap \fT', s')$ and $( \al\cap \fB_e,s' )$ are the corresponding tangle diagrams where the states on $e'$ are given by $s'$. For $e\in \{a,b,c\}$ let $\bk_{s'}(e) = \sum_{x\in \al \cap e'} s'(x)$. By the charge conservation property, Equation \eqref{eq.charge}, one of the scalars in the square bracket is 0 unless $\bk_{s'} = \bk$. From Formula \eqref{eq.triangle}, the element in the curly bracket is proportional to $y^{\bk_s}$. Hence we have \eqref{eq.tri1}.

When all the states are positive, all values of $s'$ must be positive to produce a nonzero term since it is the only state satisfying the charge conservation. In this case, none of the components of $\al\cap \fT'$ is a bad arc. Hence the element in the curly bracket is $q$-proportional to $y^{\bk}$. By Lemma \ref{r.plus} the scalar in the square bracket is a power of $q^{1/2}$. The lemma follows.
\end{proof}

\def\fBp{\fB^\bullet}
\def\cSe{\cS^\bullet}
\def\SeB{ \cSe(\fB) }
\def\Yd{\cY^\bullet}
\def\YdB{\cY^\bullet(\fB)}
\def\YdT{\cY(c)}
\def\trd{\tr_c}

By removing $a$ from $\fT$ we get a bigon, which, for the purpose of later identification, is denoted by $\fB_c$. Under the quantum trace $\tr_\fT$, the counit $\ve:\cS(\fB_c)\to \cR$ lifts to a $\cR$-linear map $\ve^\ast: \YT\to \cR$ defined by
\begin{equation}\label{eq.ves}
\ve^\ast(y^\bk)= \begin{cases}
1, &\text{if }\bk(c)= \bk(b)\text{ and } \bk(a)=0, \\
0, &\text{otherwise}.
\end{cases}
\end{equation}
That is, for $x\in \cS(\fB_c)$ we have
\begin{equation}\label{eq.ve}
\ve^\ast(\tr_\fT(x)) = \ve(x).
\end{equation}

\def\ST{\cS(\fT)}
\def\SrT{\cS^\rd(\fT)}
\def\YT{\cY(\fT)}

\subsection{Bonahon-Wong's quantum trace} \label{sec.BW}
By splitting $\fS$ along interior edges $a\in \Do$ we get a collection $\cF$ of ideal triangles. 
Every interior edge $a$ becomes two edges $a', a''$ after splitting in the set of all edges of the triangles in $\cF$.
The quantum trace map is the composition
\begin{equation}
\tr_\Delta: \CS \xhookrightarrow{\theta}
\bigotimes_{\fT\in \cF} \cS(\fT) \xhookrightarrow{\ot \tr_\fT} \bigotimes_{\fT\in \cF} \YT,
\label{eqn-tr-def}
\end{equation}
where $\theta$ is the splitting homomorphism. The image is contained in a subalgebra isomorphic to $\Ybl(\cF;\Delta)$.

Let us explain how $\Ybl(\cF;\Delta)$ embeds in the target space $A:=\bigotimes_{\fT\in \cF} \YT$. Each $\fockY(\mathfrak{T})$ is a subalgebra of $A$ by the obvious embedding. If $a\in \Do$ split into two edges $a', a''$, let $y_a:=[y_{a'} y_{a''}]$. If $a\in \Dd$, then $a$ does not split and remains an edge of a triangle in $\cF$.
Now for every $a\in \Delta$ we have an element $y_a \in A$. The subalgebra of $A$ generated by $y_a^{\pm 1}, a \in \Delta$ is isomorphic to $\cY(\fS;\Delta)$, and will be identified with $\cY(\fS;\Delta)$.

\subsection{Proof of Theorem~\ref{thm.embed3}} \label{sec.coord}

\begin{proof}
(a) Let $\fS^\ast$ be the result of attaching an ideal triangle $\fT_e$ to each boundary edge $e$ of $\fS$ by identifying $e$ with an edge of $\fT_e$. Denote the other two edges of $\fT_e$ are by $\hat e , \hat e'$ as in Figure~\ref{fig:extSurf}.
There is a smooth embedding $\iota: \fS \embed \fS^\ast$ which maps $e$ to $\hat e$ and is identity outside a small neighborhood of $e$ for every boundary edge $e$, see Figure~\ref{fig:extSurf}. In particular $\iota(a)=a$ for all $a\in \Do$.

\begin{figure}[h]
\centering
\begin{tikzpicture}
\fill[gray!20!white] (1.5,0)rectangle(-1.5,1.2) (1,0)--(-0.5,-1)--(-1,0);
\draw (-1.5,0)--(1.5,0) (-0.5,-1) edge node[below]{$\hat{e}$} (1,0)
	edge[dashed] node[left]{$\hat{e}'$} (-1,0);
\draw[fill=white] (-1,0) circle(2pt) (1,0) circle(2pt) (-0.5,-1) circle(2pt);
\draw (1.5,1.2)node[below left]{$\mathfrak{S}$};
\draw (0,0)node[above]{$e$} (-0.2,-0.4)node{$\mathfrak{T}_e$};

\begin{scope}[xshift=5cm]
\fill[gray!20!white] (1.5,0)rectangle(-1.5,1.2);
\draw (-1.5,0)--(1.5,0);
\draw[dashed] (-1,0)--(0,0);
\draw[fill=white] (-1,0) circle(2pt) (1,0) circle(2pt);
\draw[very thick] (-0.4,0.8)--(-0.4,0) (0,0.8)--(0,0) (0.4,0.8)--(0.4,0);
\draw (1.5,1.2)node[below left]{$\mathfrak{S}$} (-0.6,0)node[below]{$e$};
\draw[thick,decorate,decoration=brace] (-0.45,0.9)--(0.45,0.9);
\draw (0,1.2)node{$\alpha$};
\draw (2,0.75)node{$\xrightarrow{\iota}$};
\end{scope}

\begin{scope}[xshift=9cm]
\fill[gray!20!white] (1.5,0)rectangle(-1.5,1.2) (1,0)--(-0.5,-1)--(-1,0);
\draw (-1.5,0)--(1.5,0) (-0.5,-1) edge node[below,pos=0.8]{$\hat{e}$} (1,0)
	edge[dashed] (-1,0);
\draw[dashed] (-1,0)--(0,0);
\draw[fill=white] (-1,0) circle(2pt) (1,0) circle(2pt) (-0.5,-1) circle(2pt);
\foreach \x in {-0.4,0,0.4}
\draw[very thick] (\x,0.8)--(\x,0)arc[radius=1-\x,start angle=-180,end angle=atan(2/3)-180];
\draw[thick,decorate,decoration=brace] (-0.45,0.9)--(0.45,0.9);
\draw (0,1.2)node{$\iota(\alpha)$};
\end{scope}

\end{tikzpicture}
\caption{Left: attaching triangle $\fT_e$. Right: the embedding $\iota:\fS\embed \fS^\ast$..}\label{fig:extSurf}
\end{figure}
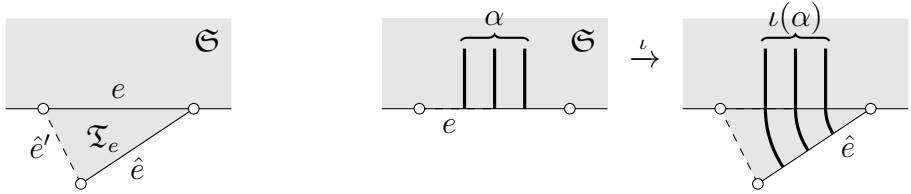


Note that $\Delta^\ast=\Delta\cup (\bigcup_{e\in \Dd} \{\hat e , \hat e'\})$ is an ideal triangulation of $\fS^\ast$. We can define the antisymmetric function $Q_{\Delta^\ast}:\Delta^\ast\times \Delta^\ast \to \BZ$, the Chekhov-Fock algebra $\cY(\fS^\ast;\Delta^\ast)$ and its balanced subalgebra $\Ybl(\fS^\ast;\Delta^\ast)$ for the surface $\surface^\ast$. We will show that the following composition
\begin{equation}\label{eq.phiD}
\phi_\Delta: \SS \xrightarrow{\iota_\ast} \cS(\fS^\ast) \xrightarrow{p} \cS^\rd(\fS^\ast) \xhookrightarrow{\tr^\rd_{\Delta^\ast}} \Ybl(\fS^\ast;\Delta^\ast)
\end{equation}
is quantum trace map we are looking for. It is important to note that under $\iota$ the image of a $\pfS$-arc is never a corner arc, and a fortiori it is not a bad arc. It follows that the composition $p\circ \iota_\ast$ maps the $\cR$-basis $B(\fS)$ of $\CS$ injectively into the basis $B^\rd(\fS^\ast)$ of $\cS^\rd(\fS^\ast)$. This shows $p\circ \iota_\ast$ is injective. Hence the composition $\phi_\Delta$ is injective.

Since $Q_{\Delta^\ast}(e, \hat e)=1=-\bQ_\Delta(e,\hat e)$, we can identify the quantum torus $\bT(\bQ)$ of \eqref{eq.TQbar} as a subalgebra of $\cY(\fS^\ast;\Delta^\ast)$ via the embedding 
\begin{equation}
\label{eq.yz}
\begin{cases}
z_a \mapsto y_a, & a\in \Do, \\
z_e \mapsto [y_e y_{\hat e}], &e\in \Dd,\\
z_{\hat e}\mapsto y_{\hat e}^{-1}, &e\in \Dd.
\end{cases}
\end{equation}

\begin{lemma}\label{r.phi1}
Suppose $\al$ is a stated simple $\pfS$-tangle diagram. With $\bbn_\al$ defined by \eqref{eq.bbn},
\begin{equation}
\label{eq.ss0}
\phi_\Delta(\al) \qeq z^{\bbn_\al} u\in\cZ(\fS;\Delta)
\end{equation}
where $u$ is a polynomial in the variables $z_a^{-2},a\in \bD$ with constant term $1$.
\end{lemma}

\begin{proof}
We can assume $\iota(\al)\cap \fS=\al$. For each boundary edge $e\in \Dd$, the intersection $\iota(\al) \cap \fT_e$ is a collection of parallel arcs with endpoints in $e$ and $\he$. For $a\in \Do$ choose an arbitrary linear order on the set $\al\cap a$, while for $e\in \Dd$, choose the order on $\iota(\al)\cap e= \al\cap e$ to be the height order of $\partial \al$. Let $\cF^\ast$ denote the set of all triangles of the triangulation $\Delta^\ast$ of $\fS^\ast$. Using the definition \eqref{eqn-tr-def} of the quantum trace via splitting, we have
\begin{equation}\label{eq.st00}
\phi_\Delta(\al)= \tr_{\Delta^\ast}(\iota(\al))= \sum_{s\in S} \prod_{\fT \in \cF^\ast}\tr_\fT(\iota(\al)\cap \fT,s),
\end{equation}
where $S$ is the set of all maps $s: \iota(\al) \cap (\bigcup_{a\in\Delta^\ast} a) \to \{\pm\}$
such that the values of $s$ on $\iota \cap \he$, $e\in \Dd$, are the states of $\iota(\al)$. The summand corresponding to $s$ in the right hand side of \eqref{eq.st00} is called the $s$-summand.
By Lemma~\ref{r.tri5}, the $s$-summand is of the form $\Coeff(s) y^{\bk_s}$, with $\Coeff(s) \in \cR$ and $\bk_s:\Delta^\ast\to \BZ$ is given by
$\bk_s(a)= \sum _{x\in \al \cap a} s(x)$ for $ a\in \Delta^\ast.$
As $s(x)=\pm 1$,
we have
\begin{equation}\label{eq.h1}
\abs{\al \cap a} - \bk_s(a) \in 2 \BN, \quad \text{ for } a\in \bD.
\end{equation}

For $e\in \Dd$, if $x\in \al\cap e$ and $s(x) <s(\iota(x))$, then the arc (in $\al\cap \fT_e$) connecting $x$ and $\iota(x)$ is a bad arc, making the $s$-summand equal to 0. Hence in \eqref{eq.st00}, $S$ can be replaced by the subset $S'$ consisting of $s$ such that $s(x) \ge s(\iota(x))$. Then for $s\in S'$,
\begin{equation}\label{eq.h2}
\bk_s(e) - \bk_s(\he) \in 2\BN \quad \text{for } e\in \Dd.
\end{equation}
By the embedding \eqref{eq.yz}, we have $y^{\bk_s} = z^{\bm_s}$, where
\begin{alignat*}{2}
\bm_s(a) &= \bk_s(a) && \text{for } a\in \Delta, \\
\bm_s(\he ) &= \bk_s(e) - \bk_s(\he)&\quad&\text{for } e\in \Dd.
\end{alignat*}
By definition \eqref{eq.bbn} of $\bbn_\al$, one has $\bbn_\al(a)= |a\cap \al|$ for $a\in \Delta$ and $\bbn_\al(\he)= |\al\cap e| - \bk_s(\he)$ for $e\in \Dd$. Hence conditions \eqref{eq.h1} and \eqref{eq.h2} imply
\begin{equation}\label{eq.z0}
\bbn_\al - \bm_s \in 2 \BN^\bD, \quad
\mathbf{m}_s(\hat{e})\ge0\text{ for } e\in \Dd.
\end{equation}
In addition, $\bm_s=\bbn_\al$ exactly when $s=s_+$, which takes value $+$ on all $a\in \Delta$. Lemma~\ref{r.tri5} shows that $\Coeff(s_+)= q^t$ for $t\in \frac 12 \BZ$. Thus, from \eqref{eq.st00} we have
\[\phi_\Delta(\al) = \sum_{s\in S'} \Coeff(s) z^{\bm_s} = q^t z^{\bbn_\al} + \sum_{s\in S'\setminus\{ s_+\}} \Coeff(s) z^{\bm_s},\]
which implies \eqref{eq.ss0} by \eqref{eq.z0} and Lemma~\ref{r.bl}.
\end{proof}

Thus
we can consider $\phi_\Delta$ as an algebra embedding
$$ \phi_\Delta: \CS\embed \cZ(\fS;\Delta).$$

When $\cR=\BZ[q^{\pm 1/2}]$, all the maps in \eqref{eq.phiD} are reflection invariant, hence so is $\phi_\Delta$.

Let us now prove that $\pr \circ \phi_\Delta=\tr_\Delta$.
Let $\cF\subset \cF^\ast$ be the subset of triangles in $\fS$, i.e., the triangles of $\Delta$. Consider the following diagram, where $\theta$ is the splitting homomorphism:
\begin{equation}\label{eq.dia2}
\begin{tikzcd}
\SS \arrow[r,"\theta'"] \arrow[d,equal]& \displaystyle{\bigotimes_{\fT\in \cF} \cS(\fT) \ot \Big(\bigotimes_{e\in \Dd } \cS(\fT_e)\Big)}
\arrow[d, "\id \ot \ve"]\arrow[r, "\ot \tr_\fT "]& \displaystyle{\bigotimes_{\fT\in \cF} \cY(\fT) \ot \Big(\bigotimes_{e\in \Dd } \cY(\fT_e) \Big)}
\arrow[d, "\id \ot \ve^\ast"]\\
\cS(\fS) \arrow[r,"\theta"] & \displaystyle{\bigotimes_{\fT\in \cF} \cS(\fT) }
\arrow[r, " \ot \tr_\fT "]& \displaystyle{\bigotimes_{\fT\in \cF} \cY(\fT)}
\end{tikzcd}
\end{equation}
Here $\ve$ is the counit of the skein algebra of the bigon $\mathfrak{B}_e=\fT_e \setminus \he'$. The left square is commutative by the counit property. The right square is commutative due to \eqref{eq.ve}. The composition of the top line is $\phi_\Delta$, while the composition of the second line is $\tr_\Delta$. Hence $(\id \ot \ve^\ast)\circ \phi_\Delta = \tr_\Delta$.

For $z^\bk\in\extYbl(\surface;\Delta)$, the definitions of $\ve^\ast$ by \eqref{eq.ves} and $\pr$ by \eqref{eq.pr} show that
$(\id \ot \ve^*)(z^\bk) = \pr(z^\bk)$.
Hence we have $\pr \circ \phi_\Delta= \tr_\Delta$. This completes the proof of part (a) of Theorem \ref{thm.embed3}.

(b) By Lemma \ref{r.phi1}, the top degree term of $\phi_\Delta(\al)$ is $q$-proportional to $z^{\bbn_\al}$, proving (b).

\def\bYt{ A^{(2)}}
\def\bYtp{ A'^{(2)}}

(c) Given another ideal triangulation $\Delta'$ of $\fS$, there is a coordinate change isomorphism
\[\Theta_{\Delta'^\ast \Delta^\ast } : \Fr(\Ybl(\fS^\ast;\Delta^\ast )) \xrightarrow{\cong} \Fr(\Ybl(\fS^\ast;\Delta'^\ast)).\]
From the compatibility with the quantum trace map, we have
\begin{equation}\label{eq.comp2}
\Theta_{\Delta'^\ast \Delta^\ast } \circ \phi_{\Delta} = \phi_{\Delta'}.
\end{equation}

Let us first prove 
\begin{equation}\label{eq.uu}
\Theta_{\Delta ^\ast \Delta'^\ast } (\Fr(\bY^\bal(\fS;\Delta))) \subset \Fr(\bY^\bal(\fS;\Delta'))
\end{equation}
Inclusion \eqref{eq.uu} follows easily from explicit formulas of $\Theta_{\Delta ^\ast \Delta'^\ast }$ given in \cite{Hiatt}, which involves many cases and are cumbersome. So we present another argument here.

Let $\bYt$ be the division subalgebra of $\Fr(\Ybl(\fS^\ast;\Delta^\ast ))$ generated by $y_a^2, a\in \bD$. Similarly let $\bYtp$ be the division subalgebra of $\Fr(\Ybl(\fS^\ast;\Delta'^\ast ))$ generated by $y_a^2, a\in \bD'$. The transition from $\Delta'$ to $\Delta$ changes only edges lying in the interior of $\fS$. Hence the locality property \cite{Bai} of the coordinate change isomorphism tells us that
\begin{align}
\Theta_{\Delta'^\ast \Delta^\ast }(y_\he^2)& = y_\he^2,\quad e\in \Dd, \label{eq.z}\\
\Theta_{\Delta'^\ast \Delta^\ast }(\bYt) &\subset \bYtp, \label{eq.zzz}
\end{align}
and the restriction of $\Theta_{\Delta'^\ast \Delta^\ast }$ on $\Fr(\bYbl(\fS;\Delta))$ is $\Theta_{\Delta' \Delta }$.

Assume $z^\bk\in \bYbl(\fS;\Delta)$. By Lemma \ref{r.bl}, there are a stated simple $\pfS$-tangle diagram $\al$ and
 $v\in \bYt$ such that
$ z^\bk= v z^{\bbn_\al}$. Formula \eqref{eq.ss0} shows that there is $v'\in \bYt$ such that
$z^{\bbn_\al}= v' \phi_{\Delta}(\al)$. Hence $$
z^\bk= vv' \phi_{\Delta}(\al).
$$
Applying $\Theta_{\Delta'^\ast \Delta^\ast}$ and use the compatibility \eqref{eq.comp2},
$$
\Theta_{\Delta'^\ast \Delta^\ast }(z^\bk)= \Theta_{\Delta'^\ast \Delta^\ast }(vv') \phi_{\Delta'}(\al).
$$
The first factor is in $\bYtp$ by \eqref{eq.zzz}, while the second factor is in $\bY^\bal(\fS;\Delta')$. This proves the inclusion \eqref{eq.uu}.

Define $\bTheta_{\Delta' \Delta }$ as the restriction of $\Theta_{\Delta'^\ast \Delta^\ast }$ to $\Fr(\bY^\bal(\fS;\Delta))$. By \eqref{eq.z} the algebra map $\bTheta_{\Delta' \Delta }$ is the extension of $\Theta_{\Delta' \Delta }$ given by $\bTheta_{\Delta' \Delta }(z_\he^2)= z_\he^2$. The compatibility \eqref{eq.comp1} follows from \eqref{eq.comp2}. Theorem \ref{thm.embed3} is proved.
\end{proof}

\subsection{Proof of Theorem~\ref{r.degen}}

\begin{proof}
We consider $\phi_\Delta:\CS \to \bT(\bQ)$ as a map of filtered algebra and consider its associated graded homomorphism. By definition the associated graded algebra is, with the convention $F_{-1}(\CS)=\{0\}$,
\[\Gr(\CS)= \bigoplus_{k=0}^\infty F_k(\CS)/F_{k-1}(\CS).\]
As $\Gr(\bT(\bQ))=\bT(\bQ)$, one has the associated algebra homomorphism
\begin{equation}\label{eq.z6}
\Gr(\phi_\Delta): \Gr(\CS) \to \bT(\bQ).
\end{equation}
Recall the set $B_k(\fS) = \{ \al \in B(\fS) \mid \dego(\al) \le k\}$ is an $\cR$-basis of $F_k(\CS)$. Hence $B_k(\fS)\setminus B_{k-1}(\fS)$ is an $\cR$-basis of $F_k(\CS)/F_{k-1}(\CS)$. It follows that $B(\fS)$ is also an $\cR$-basis of $\Gr(\CS)$.

By \eqref{eqn-tr-highest}, for any $\al\in B(\fS)$, there is $t(\al) \in \frac 12 \BZ$ such that
\begin{equation}\label{eq.z5}
\Gr(\phi_\Delta)(\al) = q^{t(\al)} z^{\bbn_\al}.
\end{equation}
By Proposition \ref{r.para},
the set $\{ z^{\bbn_\al} \mid \al\in B(\fS)\}$ is an $\cR$-basis of the monomial algebra $A(\bQ;\Lambda_\Delta)$. Equation \eqref{eq.z5} shows that $\Gr(\phi_\Delta)$ map the $\cR$-basis $B(\fS)$ of $\Gr(\CS)$ bijectively onto an $\cR$-basis of $A(\bQ;\Lambda_\Delta)$. Hence the map $\Gr(\phi_\Delta)$ of \eqref{eq.z6} is an isomorphism. This proves Theorem~\ref{r.degen}.
\end{proof}

\begin{remark}
This proof does not rely on the injectivity of $\phi_\Delta$. Thus it can serve as a proof of injectivity of $\phi_\Delta$ without using the reduced skein algebra.
\end{remark}

\subsection{Proof of Theorem~\ref{thm-noether1}}

\begin{proof}
First we assume that $\fS$ is triangulable, with an ideal triangulation $\Delta$.

(a) The submonoid $\Lambda_\Delta \subset \BN^\bD$ is finitely generated as an $\BN$-module since it is a pointed integral polyhedral cone. Let $\bk_1, \dots, \bk_n$ be a minimal set of generators of $\BN$-module $\Lambda_\Delta$. Then each $\bk_j$ cannot be the sum of two non-zero elements of $\Lambda_\Delta$.
Let $\alpha_i\in B(\fS)$ be the element such that $\bbn_{\alpha_i}=\bk_i$. It follows that each $\alpha_i$ has one component. The set $\{ (z^{\bk_1})^{m_1} (z^{\bk_2})^{m_2} \dots (z^{\bk_n})^{m_n} \mid m_i \in \BN\}$ spans $A(\bQ;\Lambda_\Delta)$ over $\BZ$. An induction on the degree using \eqref{eq.z5} shows that the corresponding set $\{ \alpha_1^{m_1} \alpha_2^{m_2} \dots \alpha_n^{m_n} \mid m_i \in \BN\}$ spans $\CS$ over $\cR$. This proves $\{\alpha_1, \dots, \alpha_n\}$ is an orderly generating set for $\stateS(\surface)$.

(b) By Lemma \ref{r.mono} the ring $A(\bQ,\Lambda_\Delta)$ is a Noetherian domain. As an associated algebra of $\CS$ is a Noetherian domain, so it $\CS$.

(c) By Lemma~\ref{lemma-GKdim}, the monomial algebra $A(\bQ;\Lambda_\Delta)$ has GK dimension $r(\fS)=\abs{\bD}$, the rank of $\Lambda_\Delta$. For each $k$, the $\cR$-module $F_k(\SS)$ is free of finite rank. Then the GK dimension of $\CS$ is the same as the GK dimension of the associated graded algebra $A(\bar{Q};\Lambda_\Delta)$, which is $r(\fS)$.

This completes the proof of the theorem for the case of triangulable surfaces. Suppose $\fS$ is not triangulable. We already know $\CS$ is a domain. By removing enough points from $\fS$, we can obtain a new surface $\fS'$ which is triangulable. As $\CS$ is a quotient of $\cS(\fS')$, we get the orderly finite generation and the Noetherian property for $\CS$.
\end{proof}

\section{Relation between the two quantum trace maps}\label{sec.rel}

\subsection{Comparing two quantum trace maps}
Assume that $\fS$ is a triangulable punctured bordered surface with non-empty boundary, with $\Po$ the set of interior punctures.

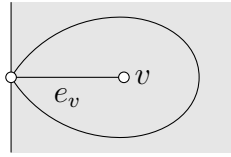
\begin{figure}[h]
\centering
\begin{tikzpicture}
\fill[gray!20!white] (0,0)rectangle(3,2);
\draw (0,1) to node[below]{$e_v$} (1.5,1)node[right]{$v$};
\draw (0,0)--(0,2) (0,1) to[curve through={(1.5,1.8)..(2.5,1)..(1.5,0.2)}] cycle;
\draw[fill=white] (0,1)circle(2pt) (1.5,1)circle(2pt);
\end{tikzpicture}
\caption{The interior puncture $v$ is in a monogon, and $e_v$ connects $v$ and the vertex of the monogon.}\label{fig:monarc}
\end{figure}

Let $\cE$ be a quasitriangulation of $\fS$. Then $\cE$
can be uniquely \emph{completed} to an ideal triangulation $\Delta=\cE \cup (\bigcup_{v\in\Po} e_v)$, where $e_v$ is an ideal arc disjoint from $\cE$, and connecting $v$ and the vertex of the monogon containing $v$, see Figure \ref{fig:monarc}.

By Theorems \ref{thm.embed1} and \ref{thm.embed3}, we have two algebra embeddings, the length coordinate quantum trace map and the shear coordinate quantum trace map
\begin{align}
\vpE&: \CS \embed \bXE, \label{eq.ll}\\
\phi_\Delta &: \CS \embed \bYbl(\fS;\Delta).
\end{align}
Here $\bXE$ is the quantum torus over $\RV$ generated by $x_e, e\in \bE$, with relations $x_a x_b = q^{\bP(a,b)} b_b x_a$ for $a,b,\in \bE$, while $\bYbl(\fS;\Delta)$ is the balanced subalgebra of $\bT(\bQ_\Delta)$.



The polynomial ring $\RV$ is considered as a subring of the center of $\CS$ by identifying $v\in \Po$ with a small loop $X_v\subset \fS$ surrounding $v$. The map $\vpE$ in \eqref{eq.ll} is an $\RV$-algebra homomorphism. Identify $\RV$ as a subalgebra of the Laurent polynomial ring $\RVb:=\cR[x_v^{\pm1}, v\in \Po]$ by $v= x_v + x_v^{-1}$. By extending the ground ring from $\RV$ to $\RVb$ from \eqref{eq.ll} we get
\begin{align}
\vpEb: \CSb \embed \bXEb, \label{eq.ll2}
\end{align}
where $\CSb= \CS\ot_\RV \RVb$ and $ \bXEb = \bXE \ot_\RV \RVb$.

\def\bPb{{\bP^\diamond }}

Note that $\bXEb$ is the quantum torus over $\cR$ associated to the antisymmetric function $\bPb: (\bar{\quasi} \cup \Po)\times (\bar{\quasi} \cup \Po) \to \BZ$ which is the 0 extension of $\bP:\bar{\quasi} \times \bar{\quasi} \to \BZ$.

We will consider $\bXE$ as a subalgebra of $\bXEb$ by the obvious embedding, and
also use $\vpE$ to denote the composition $\CS \xhookrightarrow{ \vpE} \bXE \embed \bXEb$.


\def\psE{\psi_\cE}
\begin{theorem}\label{thm-len2shear}
Let $\fS$ be a triangulable connected punctured bordered surface with non-empty boundary, $\cE$ be a quasitriangulation of $\fS$, and $\Delta$ be the unique triangulation containing $\cE$.

There is a unique multiplicatively linear 
$\cR$-algebra isomorphism
\[\psi_\quasi:\bXEb\xrightarrow{\cong} \extYbl(\surface;\Delta)\]
such that $\psE(x_v) = y_{e_v}$ for $v\in \Po$ and the following diagram commutes:
\[\begin{tikzcd}[row sep=tiny]
&\extYbl(\surface;\Delta)\\
\stateS(\surface)\arrow[ru,"\phi_\Delta"]
\arrow[rd,"\varphi_\quasi"']&\\
&\extX^\diamond (\surface;\quasi)\arrow[uu,"\psi_\quasi"]
\end{tikzcd}\]
\end{theorem}
Both $\psE$ and its inverse have transparent, simple geometric interpretation, which is explained in the proof.

\def\pDb{\pD^\diamond }

\begin{corollary}\label{r.iso5}
Suppose $\fS$ is a triangulable connected punctured bordered surface with non-empty boundary, and $\Delta$ is an ideal triangulation of $\fS$, not necessarily coming from a quasitriangulation. The shear coordinate quantum trace $\pD: \CS \embed \bYD$ can be extended to an isomorphism of division algebras
\begin{equation}
\tilde{\phi}_\Delta: \Fr(\CSb) \xrightarrow{\cong} \Fr(\bYD)
\end{equation}
such that the coordinate change isomorphism is given by $\bTheta_{\Delta \Delta'}= \tilde{\phi}_\Delta\circ\tilde{\phi}_{\Delta'}^{-1}$.
\end{corollary}

\begin{proof}
Suppose $\Delta_1$ is the completion of a quasitriangulation $\cE$. Define $\tilde{\phi}_{\Delta_1}$ as the composition
\[\Fr(\CSb) \xrightarrow{\Fr(\vpEb)} \Fr(\bXEb) \xrightarrow{\Fr(\psE)} \Fr(\extYbl(\surface;\Delta_1)),\]
which is an isomorphism by Theorem~\ref{thm-len2shear}. For any triangulation $\Delta$ define $\tilde{\phi}_{\Delta}= \bTheta_{\Delta\Delta_1}\circ\tilde{\phi}_{\Delta_1}$. It is trivial to verify the properties of $\tilde{\phi}_{\Delta}$.
\end{proof}

If in addition, $\fS$ has no interior punctures, then $\bE=\bD$ and $\bpXE=\bXE$.
\begin{corollary}\label{r.iso6}
Suppose $\fS$ is a triangulable connected punctured bordered surface with non-empty boundary and no interior punctures, and $\Delta$ is an ideal triangulation of $\fS$. The algebra embedding $\pD: \CS \embed \bYD$ induces an isomorphism $\Fr(\pD): \Fr (\CS) \xrightarrow{\cong} \Fr \bYD$ of division algebras.
\end{corollary}
With the assumption of Corollary~\ref{r.iso6}, one can prove the existence of the coordinate change isomorphism, which is much simpler than the original construction: If $\Delta'$ is another ideal triangulation, then define $\bar{\Theta}_{\Delta\Delta'}= \Fr(\phi_{\Delta}) \circ \Fr(\phi_{\Delta'})^{-1}$.

\begin{remark}
If there are interior punctures and potentially empty boundary, we can use the ``open up" procedure in \cite{FT,Le:qtrace}, which can give an alternative definition of $\bar{\Theta}_{\Delta\Delta'}$. The details require a significant detour, so we are not going to explore this possibility here.
\end{remark}

\subsection{Proof of Theorem \ref{thm-len2shear}}

Recall that for $a\in \bE$ there is a normalized stated arc $X_a$ mapping to $x_a$ under $\vpE$, see Subsection \ref{sec.bad}. It is not hard to show that $\pD(X_a)$ is a monomial as in the following.

\begin{lemma}\label{lemma-matrix-K}
(Compare Lemma~\ref{r.phi1})
For $a\in \bE$ one has $ \pD(X_a)= y^{K_a}$, where $K_a:\bD \to \BN$ is given by
\begin{equation}\label{eq.Darc}
K_{a}(c)=\#\left( \vertexdef{$c$}{$a$} \right)+
\begin{cases}
-2,&\text{if }a=\hat{e}=c\text{ for some }e\in \Dd,\\
0,&\text{otherwise}.
\end{cases}
\end{equation}
\end{lemma}

\begin{remark}
The picture in \eqref{eq.Darc} appeared in Subsection~\ref{sec.P.arc} when both arcs are boundary ending. The definition still makes sense as long as the arcs do not meet at an interior puncture. Since $a$ is boundary ending, $K_a(c)$ is well defined. Note when $a=c,\hat{c}$ or vice versa, the correct count requires an isotopy so that $a$ and $c$ are disjoint. 
\end{remark}

\begin{proof}
Recall the element $X_a$ is defined by assigning states to the arc $D(a)$, which is obtained by moving the endpoints of $a$ to the left. We also need to apply the map $\iota:\surface\to\surface^\ast$ to compute $\phi_\Delta(X_a)$.
To draw the arc $\iota(X_a)$, near each end of the ideal arc $a$, start slightly to the left of $a$ and draw a segment in the counterclockwise direction all the way to the boundary of $\surface^\ast$. Now connect the two segments, which crosses $a$ at one point since the two segments are on different sides of $a$. This is shown in Figure~\ref{fig-arc-matrix}. (Technically the corners should also be smoothed.)

\begin{figure}[h]
\centering
\begin{subfigure}[b]{0.3\linewidth}
\centering
\begin{tikzpicture}
\coordinate (O) at (0,0);
\coordinate (E) at (90:2);
\fill[gray!20!white] (1.2,0) rectangle (-1.2,2)
	(O) -- ++(-150:1.2) -- (-1.2,0) (E) -- ++(30:1.2) -- (1.2,2);
\draw (O) edge node[below]{$\hat{e}_1$} ++(-150:1.2)
	edge node[above,pos=0.8]{$e_1$} (180:1.2) edge (150:0.6) edge (120:0.6)
	edge (E) edge (30:0.6) edge (0:1.2);
\draw (E) edge node[above]{$\hat{e}_2$} ++(30:1.2)
	edge node[below,pos=0.8]{$e_2$} ++(0:1.2) edge ++(-30:0.6) edge ++(-60:0.6)
	edge ++(-150:0.6) edge ++(180:1.2);
\draw (0,1) node[left]{$\iota(X_a)$} (0,0.7) node[right]{$a$};
\draw[inner sep=0pt] (-150:0.4)node[below right]{$+$}
	(E) ++ (30:0.4)node[above left]{$+$};
\draw[very thick] (E) ++(30:0.4)
	arc[radius=0.4,start angle=30, end angle=-75]
	--(105:0.4)
	arc[radius=0.4,start angle=105, end angle=210];
\draw[fill=white] (O)circle(2pt) (E)circle(2pt);
\end{tikzpicture}
\subcaption{$a\in\mathring{\Delta}$}
\end{subfigure}
\begin{subfigure}[b]{0.3\linewidth}
\centering
\begin{tikzpicture}
\fill[gray!20!white] (1.5,0) -- ++(-135:1.2)--(0,0)--(-135:1.2)--(-1.5,0)
	(1.5,0)arc[start angle=0,end angle=180,radius=1.5];
\draw[pos=0.7] (-1.5,0)--(1.5,0)
	(1.5,0)node[above left]{$e$} (75:0.5)node[above]{$\iota(X_a)$}
	(0,0) edge node[left]{$\hat{e}_1$} (-135:1.2)
	edge (30:0.7) edge (120:0.7) edge (150:0.7)
	edge node[above]{$e_1$}(-1.5,0)
	(1.5,0) edge node[left]{$\hat{e}$} ++(-135:1.2);
\draw[very thick,inner sep=0pt] (1.5,0)++(-135:0.4)node[below right]{$\pm$}
	arc[radius=0.4,start angle=-135,end angle=-165]
	--(15:0.4)arc[radius=0.4,start angle=15,end angle=225]node[below right]{$+$};
\draw[fill=white] (0,0) circle(2pt) (1.5,0)circle(2pt);
\end{tikzpicture}
\subcaption{$a=e,\hat{e}$ for $e\in\mathring{\Delta}$}
\end{subfigure}
\caption{The arc $\iota(X_a)$.}\label{fig-arc-matrix}
\end{figure}
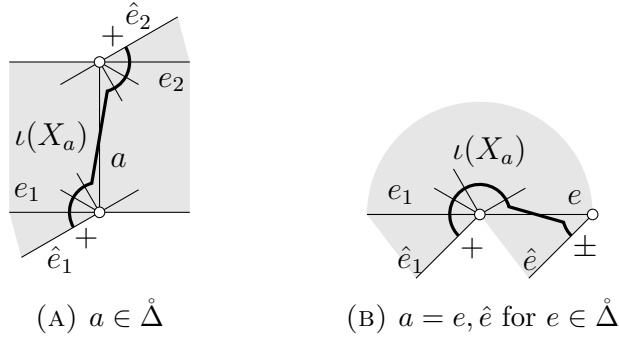

Fix a direction of $X_a$ and choose the height to be increasing following the direction. In the state sum formula \eqref{eq.st00}, where $\al= X_a$, the only $s$ which does not result in a bad arc is the one taking only value $+$ on $X_a\cap (\bigcup_{b\in \Delta} b)$. This can be seen by starting at the endpoint with $+$ state. Each time the arc crosses an edge of $\Delta$ up to the point where it meets $a$, avoiding bad arcs forces the sign to be $+$. This only term gives a monomial in $y$, and by the description above, the exponent is given by $K_a$ of \eqref{eq.Darc}. 
Thus, $\pD(X_a)\qeq y^{K_a}$. Then reflection invariance forces $ \pD(X_a)= y^{K_a}$.
\end{proof}

For the generators $x_a$ of the quantum torus $\bpXE$ define $\psE(x_a) = \pD(X_a)$ for $a\in \bE$, and define $\psE(x_v)= y_{e_v}$ for $v\in \Po$. We have $\psE(x_a x_b)= \psE(x_a) \psE(x_b)$ for $a,b\in \bE$ because $\pD$ is an algebra homomorphism. The remaining variables $x_v, v\in \Po$, are central, and their images $y_{e_v}$ are also central. Thus all the defining relations among the variables of the quantum torus $\extYbl(\surface;\Delta)$ are respected, and we can extend $\psE$ to an algebra homomorphism $\psE: \bpXE \to \bYD $, which is multiplicatively linear because the image of each variable is a monomial. The restriction of $\psE$ to $\bXE$ is $\pD$, since $\pD(v)= y_{e_v} + y_{e_v}^{-1}$ by the state sum formula \eqref{eq.st00} with $\al= X_v$.

The map $\psE$ is injective since otherwise its image is a quantum torus of fewer variables, i.e., of GK dimension less than $r(\fS)$, but the image of $\psE$ contains the image of the injective map $\pD$, which has GK dimension $r(\fS)$.

Let us now prove that $\psE$ is surjective. Recall by Lemma~\ref{r.bl}, the algebra $\bYD$ is generated by elements of the form $z^{\bar{\mathbf{n}}_\alpha}$ with $\al$ a stated simple $\pfS$-tangle diagram and $z_e^{\pm2}$ with $e\in \bD$. As $\al\in \stateS(\surface)\subset \bXE$, $\pD(\al)=\psE(\al)$ is in the image of $\psE$. Since the image is a monomial subalgebra, the top degree monomial, $z^{\bar{\mathbf{n}}_\alpha}$ by Lemma~\ref{r.phi1}, is also in the image. It remains to show that $z_e^2$, or alternatively $y_e^2$, is in the image of $\psE$. Define $\bar{H}_e:\bar{\quasi}\cup\mathring{\marked}\to\mathbb{Z}$ for each $e\in\bar{\Delta}$ by
\begin{itemize}
\item $x^{\bar{H}_{e_v}}:=x_v^2$ for $v\in\mathring{\marked}$,
\item $x^{\bar{H}_e}:=[x_a x_b^{-1} x_c x_d^{-1}]$ for $e\in \Do \setminus \Po$, where $e,a,b$ and $e,c,d$ are triangular triples in counterclockwise order,
\item 
$x^{\bar{H}_e}:= [x_a x_b^{-1} x_{\hat{e}}]$, $x^{\bar{H}_{\hat{e}}}:=[x_e x_\he^{-1}]$ for $e \in \Dd$, where $e,a,b$ form a triangular triple in counterclockwise order.
\end{itemize}
We will show
\begin{equation}\label{eq.zt1}
\psE(x^{\bar{H}_e})=y^{\sigma_e}:=\begin{cases}
y_e^2y_{e_v},&e\text{ bounds a monogon containing }v\in\mathring{\marked},\\
y_e^2,&\text{otherwise}.
\end{cases}
\end{equation}
This is trivial when $e\in\mathring{\marked}$. For $e\in\bar{\quasi}$, $\bar{H}_e(\mathring{\marked})=0$. Thus \eqref{eq.zt1} is equivalent to
\begin{equation}\label{eq.eq1}
\sigma_e(c) = \sum_{a\in \bar{\quasi}} \bar{H}_e(a) K_a(c),\quad c\in\bar{\Delta},
\end{equation}
which is proved after Lemma~\ref{lemma-H-inverse} in the Appendix. This completes the proof of Theorem~\ref{thm-len2shear}.

\appendix
\section{Triangulation and quasitriangulation}

\subsection{Face and vertex matrices}

Here we give a more formal definition of the face and vertex matrices.

Suppose $\Delta$ is an ideal triangulation of the punctured bordered surface $\surface$. For each non-self-folded face $\mathfrak{T}\in\face$, the edges of $\mathfrak{T}$ are cyclically ordered. Given $a,b\in\Delta$, if $a$ and $b$ are distinct edges of $\mathfrak{T}$, define
\[Q_\mathfrak{T}(a,b)=\begin{cases}
1,&b\text{ is clockwise to }a,\\
-1,&b\text{ is counterclockwise }a.\\
\end{cases}\]
Otherwise, let $Q_\mathfrak{T}(a,b)=0$. For a self-folded face $\mathfrak{T}$, set $Q_\mathfrak{T}=0$. Then the face matrix is
\begin{equation}\label{eqn-Q-sum}
Q=\sum_{\mathfrak{T}\in\face}Q_\mathfrak{T}.
\end{equation}

The vertex matrix is defined for a quasitriangulation $\quasi$.
For each boundary puncture $v\in\marked_\partial$ and half-edges $a',b'$, let $P_{+,v}(a',b')=1$ if $a'\ne b'$ and $b'$ is counterclockwise to $a'$ at $v$. Otherwise let $P_{+,v}(a',b')=0$. Define the ``positive part" $P_+$ of the vertex matrix by
\[P_+(a,b)=\sum P_{+,v}(a',b'),\]
where the sum is over punctures $v\in\mathring{\marked}$ and half-edges $a',b'$ of arcs $a,b\in\quasi$. For this definition to make sense, an isotopy might be necessary to make $a$ and $b$ disjoint, which is important for the diagonal elements. The vertex matrix is then $P=P_+-P_+^T$. Alternatively, let $P_v(a',b')=P_{+,v}(a',b')-P_{+,v}(b',a')$. Then $P(a,b)=\sum P_v(a',b')$.

\subsection{Extended matrices}

Define the projection matrix $J:\Delta\times\Delta_\partial\to\mathbb{Z}$ or $J:\quasi\times\quasi_\partial\to\mathbb{Z}$ (depending on the context) by
\[J(a,b)=\begin{cases}1,&a=b,\\0,&a\ne b.\end{cases}\]

The extended matrices $\bar{Q}:\bar{\Delta}\times\bar{\Delta}\to\mathbb{Z}$ and $\bar{P},\bar{P}_+:\bar{\quasi}\times\bar{\quasi}\to\mathbb{Z}$ can be written in block matrix form
\begin{align}
\bar{Q}&=\begin{pmatrix}Q&-J\\J^T&0\end{pmatrix},\quad
\bar{Q}^\ast=\begin{pmatrix}Q&J\\-J^T&0\end{pmatrix}\\
\bar{P}&=\begin{pmatrix}P&-(P_++P_+^T)J\\J^T(P_++P_+^T)&-J^T PJ\end{pmatrix},\\
\bar{P}_+&=\begin{pmatrix}P_+&P_+J\\J^T P_+&J^T P_+J-2I\end{pmatrix}.\label{eqn-Pp-def}
\end{align}
Since $J$ is a projection, it picks out blocks of matrices with the correct dimension. Here we also included the restriction $\bar{Q}^\ast:\bar{\Delta}\times\bar{\Delta}\to\mathbb{Z}$ of $Q_{\Delta^\ast}$, which can be used to define $\extYbl(\surface;\Delta)$ using $y$ coordinates.

\subsection{Face-vertex matrix duality}

For a punctured bordered surface $\surface$ with at least one boundary puncture, there is a relation between the face matrix and the vertex matrix. Consider a quasitriangulation $\quasi$ and its completion $\Delta$. Let $H=I_\partial-Q_\quasi$, where $Q_\quasi$ is the restriction to $\quasi\times\quasi$, and $I_\partial=JJ^T$. The restriction to $\quasi$ does not lose any information, since by definition, all entries involving a self-folded edge is $0$, so the full matrix can be recovered by a $0$-extension, i.e., $Q=Q_\quasi\oplus 0$. The extended matrix $\bar{H}$ is defined by, using the block matrix form,
\[\bar{H}=\begin{pmatrix}-Q_\quasi&J\\J^T&-I\end{pmatrix}.\]

\begin{lemma}\label{lemma-H-inverse}
$HP_+=2I$, $\bar{H}\bar{P}_+=2I$. 
\end{lemma}

\begin{proof}
The extended case follows from the first part by a routine calculation.

The punctured monogon is a special case we can directly verify. Here $P_+=(2)$ and $H=(1)$. Thus we focus on other surfaces from now on. In particular, all edges that bound self-folded faces are interior.

There are three cases to consider.

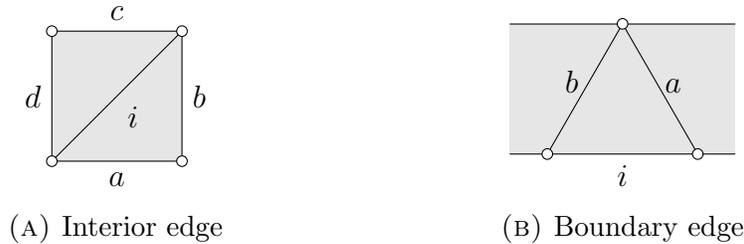
\begin{figure}[h]
\centering
\begin{subfigure}[b]{0.4\linewidth}
\centering
\begin{tikzpicture}
\coordinate (A) at (0,0);
\coordinate (B) at (1.73,0);
\coordinate (C) at (1.73,1.73);
\coordinate (D) at (0,1.73);
\draw[fill=gray!20!white] (C) to node[below right]{$i$} (A) to node[below]{$a$} (B) to node[right]{$b$} (C) to node[above]{$c$} (D) to node[left]{$d$} (A);
\draw[fill=white] (A)circle(2pt) (B)circle(2pt) (C)circle(2pt) (D)circle(2pt);
\end{tikzpicture}
\subcaption{Interior edge}\label{fig-quad-i}
\end{subfigure}
\begin{subfigure}[b]{0.4\linewidth}
\centering
\begin{tikzpicture}
\coordinate (A) at (0:2);
\coordinate (B) at (60:2);
\coordinate (C) at (0,0);
\coordinate (D) at (2.5,1.73);
\coordinate (E) at (-0.5,1.73);
\coordinate (F) at (-0.5,0);
\coordinate (G) at (2.5,0);
\fill[gray!20!white] (D)--(E)--(F)--(G);
\draw (A)--(B)--(C);
\draw (1,0)node[below]{$i$} (30:1.93)node{\vphantom{$b$}$a$}
	(60:1)++(150:0.2)node{$b$};
\draw (F)--(G) (D)--(E);
\draw[fill=white] (A)circle(2pt) (B)circle(2pt) (C)circle(2pt);
\end{tikzpicture}
\subcaption{Boundary edge}\label{fig-tri-i}
\end{subfigure}
\caption{Local pictures of the computation}
\end{figure}

\textbf{Case 1}. Suppose $i\in\mathring{\quasi}$ does not bound a self-folded face. By definition,
\[(HP_+)(i,j)=\sum_{k\in\quasi}H(i,k)P_+(k,j)=P_+(a,j)-P_+(b,j)+P_+(c,j)-P_+(d,j),\]
where $a,b,c,d$ are as in Figure~\ref{fig-quad-i}. Some of the sides might coincide, but using the sum-of-faces definition \eqref{eqn-Q-sum} of $Q$, the result is unaffected. We further split the sum into half-edges, and group the half-edges into pairs by corners of the quadrilateral.

First we consider $j\notin\{i,a,b,c,d\}$. If a half-edge $j'$ ends on one of the boundary punctures of the quadrilateral, then every corner incident to that puncture is entirely on one side of $j'$. Thus every such corner contributes $1-1$ or $0-0$ to the sum, so the total is $0$.

If $j\in\{a,b,c,d\}$, isotope $j$ into the interior of the quadrilateral. $j$ goes through two consecutive corners, and they contribute $1-1$ to the sum. In the two corners that $j$ does not go through, the calculation is as before and they cancel out.


If $j=i$, then $a$ and $c$ both contribute $+1$ in the local picture, and any other corner incident to either end of $i$ contributes $0$ as before. Thus the total is $2$.


\textbf{Case 2}. Suppose $i\in\mathring{\quasi}$ bounds a self-folded face. We can use the same picture Figure~\ref{fig-quad-i} and identify $c=d$. In this case
\[(HP_+)(i,j)=P_+(a,j)-P_+(b,j).\]
$a$ and $b$ form a (punctured) bigon. We can group half-edges of $a$ and $b$ into corners and repeat the calculation as in Case 1.

\textbf{Case 3}. Suppose $i\in\quasi_\partial$, then
\[(HP_+)(i,j)=P_+(i,j)+P_+(a,j)-P_+(b,j),\]
where $a,b$ are as in Figure~\ref{fig-tri-i}. The calculation is again analogous to the first case. Note even though the $i$ and $a$ terms have the same sign, there is nothing counterclockwise to the $i,a$ corner.
\end{proof}

\begin{proof}[Proof of \eqref{eq.zt1}]
This is a matrix equation, where $\sigma,K:\bar{\quasi}\times\bar{\Delta}\to\mathbb{Z}$ are given by
\[\sigma(a,c)=\sigma_a(c),\quad K(a,c)=K_a(c).\]
Then \eqref{eq.zt1} can be written as $\sigma=\bar{H}K$. By Lemma~\ref{lemma-H-inverse}, we just need to show $\bar{P}_+\sigma=2K$, i.e.
\begin{equation}\label{eqn-P-K}
\sum_{b\in\bar{\quasi}}\bar{P}_+(a,b)\sigma_b(c)=2K_a(c),\quad a\in\bar{\quasi},c\in\bar{\Delta}.
\end{equation}
Also note $\bar{P}_+$ is the restriction of $K$ to $\bar{\quasi}\times\bar{\quasi}$ by definition. Compare \eqref{eq.Darc} with \eqref{eqn-Pp-def}.

If $c\in\bar{\quasi}$, then $\sigma_b(c)$ is nonzero only when $b=c$, and $\sigma_c(c)=2$. Thus the equation is satisfied.

If $c=e_v$ for $v\in\mathring{\marked}$, then $\sigma_b(e_v)$ is nonzero only when $b=b_v$ is the arc bounding the monogon containing $v$. In this case, $\sigma_{b_v}(e_v)=1$. \eqref{eqn-P-K} is then equivalent to
\[\bar{P}_+(a,b_v)=2K_a(e_v).\]
This is true by checking the local picture at $v$. We can always isotope $a$ to be outside of the monogon. Then for each half edge of $a$, either $b_v$ counts twice and $e_v$ counts once, or none of them counts. Thus we have the desired equality.
\end{proof}

\end{document}